\documentclass[preprints,article,accept,moreauthors,pdftex]{Definitions/mdpi}

\usepackage{verbatim}
\usepackage{amsmath,amssymb,amsthm}             %数学公式
\usepackage{amssymb}             %数学公式
\usepackage{amsfonts}            %数学字体
\usepackage{mathrsfs}            %数学花体 
\usepackage{amsthm}
\usepackage{algorithm}  
\usepackage{algorithmicx}  
\usepackage{algpseudocode}  
\usepackage{mathtools}
\usepackage{commath}
\usepackage{physics}
\usepackage{bm}
\usepackage{graphicx}
\usepackage{geometry}
\usepackage{float}
\usepackage{listings}
\usepackage{subfigure}
\usepackage{multirow}
\usepackage{diagbox}
\usepackage{color}
\usepackage[export]{adjustbox}
\usepackage{ulem}
\usepackage{geometry}                % See geometry.pdf to learn the layout options. There are lots.
\usepackage{graphicx}
\usepackage{amssymb}
\usepackage{color}
\usepackage{bbm, dsfont}
\usepackage{latexsym,amssymb,amsmath,amsfonts,amsthm}
\usepackage{amsthm,amsmath}
\usepackage{float}
\usepackage{epic}
\usepackage{fullpage}
\usepackage{bbm}
\usepackage{subfigure}
\usepackage{graphicx}
\usepackage[toc,page]{appendix}
\usepackage{amsfonts}
\usepackage[all]{xy}
\usepackage{caption}
\usepackage{geometry}              
\usepackage{graphicx}
\usepackage{amssymb}
\usepackage{color}
\usepackage{bbm, dsfont}
\usepackage[utf8]{inputenc}
\usepackage{tikz}
\usetikzlibrary{decorations.pathmorphing}
\tikzset{snake it/.style={decorate, decoration=snake}}
\usetikzlibrary{shapes.geometric,positioning,decorations.pathreplacing}

\usepackage{amsmath}
\usepackage{mathrsfs}

\newcommand{\BZ}{{\mathbb{Z}}}

\numberwithin{equation}{section}

\newcommand{\ind}{{\mathbbm{1}}}

\newtheorem{question}{Question}

\newtheorem{conjecture}[theorem]{Conjecture}
\newcommand{\RNum}[1]{\uppercase\expandafter{\romannumeral #1\relax}}

\firstpage{1} 
\pubvolume{xx}
\issuenum{1}
\articlenumber{5}
\pubyear{2020}
\copyrightyear{2020}
\history{Received: date; Accepted: date; Published: date}
\Title{On (non-)monotonicity and phase diagram of finitary random interlacement}

\Author{Zhenhao Cai $^{1,\ddagger}$, Yunfeng Xiong $^{1,\ddagger}$ and Yuan Zhang $^{3,*} \orcidA{}$}

\AuthorNames{Zhenhao Cai, Yunfeng Xiong and Yuan Zhang}

\address{%
$^{1}$ \quad Peking University; caizhenhao@pku.edu.cn\\
$^{2}$ \quad Peking University; xiongyf@math.pku.edu.cn\\
$^{3}$ \quad Peking University; zhangyuan@math.pku.edu.cn}

\corres{Correspondence: zhangyuan@math.pku.edu.cn}

\secondnote{These authors contributed equally to this work.}
\abstract{In this paper, we study the evolution of a Finitary Random Interlacement (FRI) with respect to the expected length of each fiber. In contrast to the previously proved phase transition between sufficiently large and small fiber length, we show that for $d=3,4$, FRI is NOT stochastically monotone as fiber length increasing. At the same time, numerical evidences still strongly support the existence of a unique and sharp phase transition on the existence of a unique infinite cluster, while the critical value for phase transition is estimated to be an inversely proportional function with respect to the system intensity. }

  % List three to ten pertinent keywords specific to the article, yet reasonably common within the subject discipline.

\begin{document}
%%%%%%%%%%%%%%%%%%%%%%%%%%%%%%%%%%%%%%%%%%

%%%%%%%%%%%%%%%%%%%%%%%%%%%%%%%%%%%%%%%%%%

\section{Introduction}

Phase transition, which qualitatively characterizes the change in the state of a system under a continuous change in an external parameter, is ubiquitously found in probabilistic models and statistical mechanics. In this paper, we investigate the phase transitions in the Finitary Random Interlacement (FRI) introduced by Bowen in his study on Gaboriau-Lyons problem \cite{bowen2019finitary}. In contrast to its profound connection with the von Neumann-Day problem, a relatively simple description of FRI is given by Bowen in \cite{bowen2019finitary} as follows: Consider a  random network $(G, V)$ in $\BZ^d, \ d\ge 3$. For each vertex $x \in V$, there lives $N_x$ frogs, where $N_x$ is a Poisson random variable with mean $u \textup{deg}_x/(T+1)$, $\textup{deg}_x$ is the degree of $x$ and $u, T$ are two positive parameters. Each frog has a coin that lands on head with probability $T/(T+1)$. At time $t = 0$, the frog flips the coin. If it lands on heads, the frog moves to a random neighboring vertex with equal probability. It repeats this operation until the coin lands on tails at which point the frog stops forever. The FRI is the random multiset of random walk paths of all frogs.

Since each path consists of a simple random walk for $t$ steps and a geometric random variable with mean $T+1$ at $t+1$ steps, a FRI can be roughly treated as a random network $(G, V)$ in $\BZ^d$ ``interlaced" by fibers made of geometrically truncated simple random walk (SRW) trajectories, with a multiplicative parameter $u$ controlling its Poisson intensity, and truncation parameter $T$ that determines the expected length of each fiber. As pointed out by an anonymous referee (of a previous paper), an FRI can also be described as a variant of the Random Interlacement (RI) \cite{Sznitman2009Vacant} in $\BZ^d$ with weight \cite{sznitman2012topics}, determined by capacity with a discrete killing measure \cite{teixeira2009interlacement}. See Section \ref{sec_def} for more precise definitions and constructions for FRI.

\begin{figure}[!h]
\centering
\subfigure[T = 1.4]{\includegraphics[width=0.49\textwidth,height=0.30\textwidth]{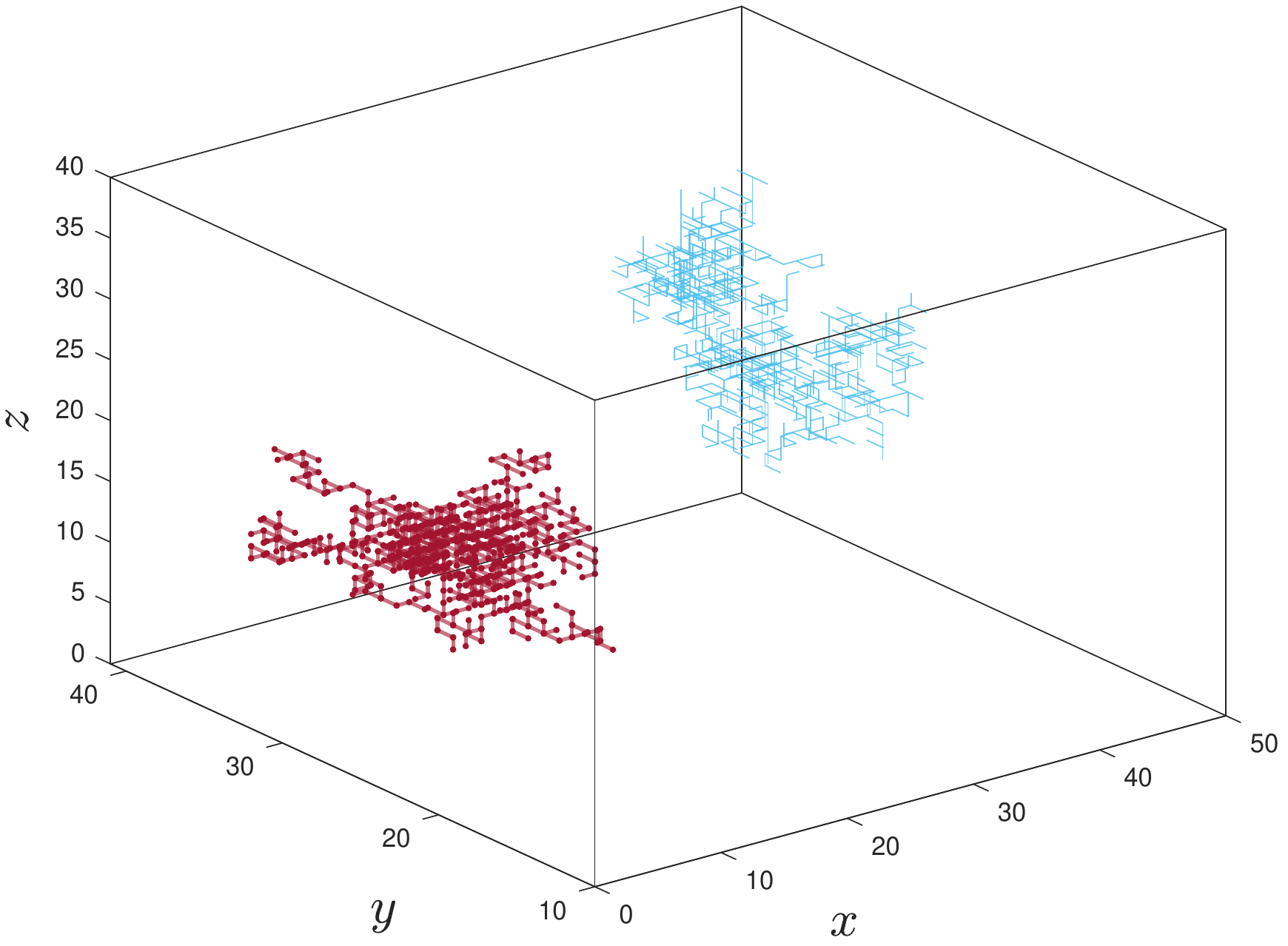}\label{fig_1a}}
\subfigure[T = 1.8]{\includegraphics[width=0.49\textwidth,height=0.30\textwidth]{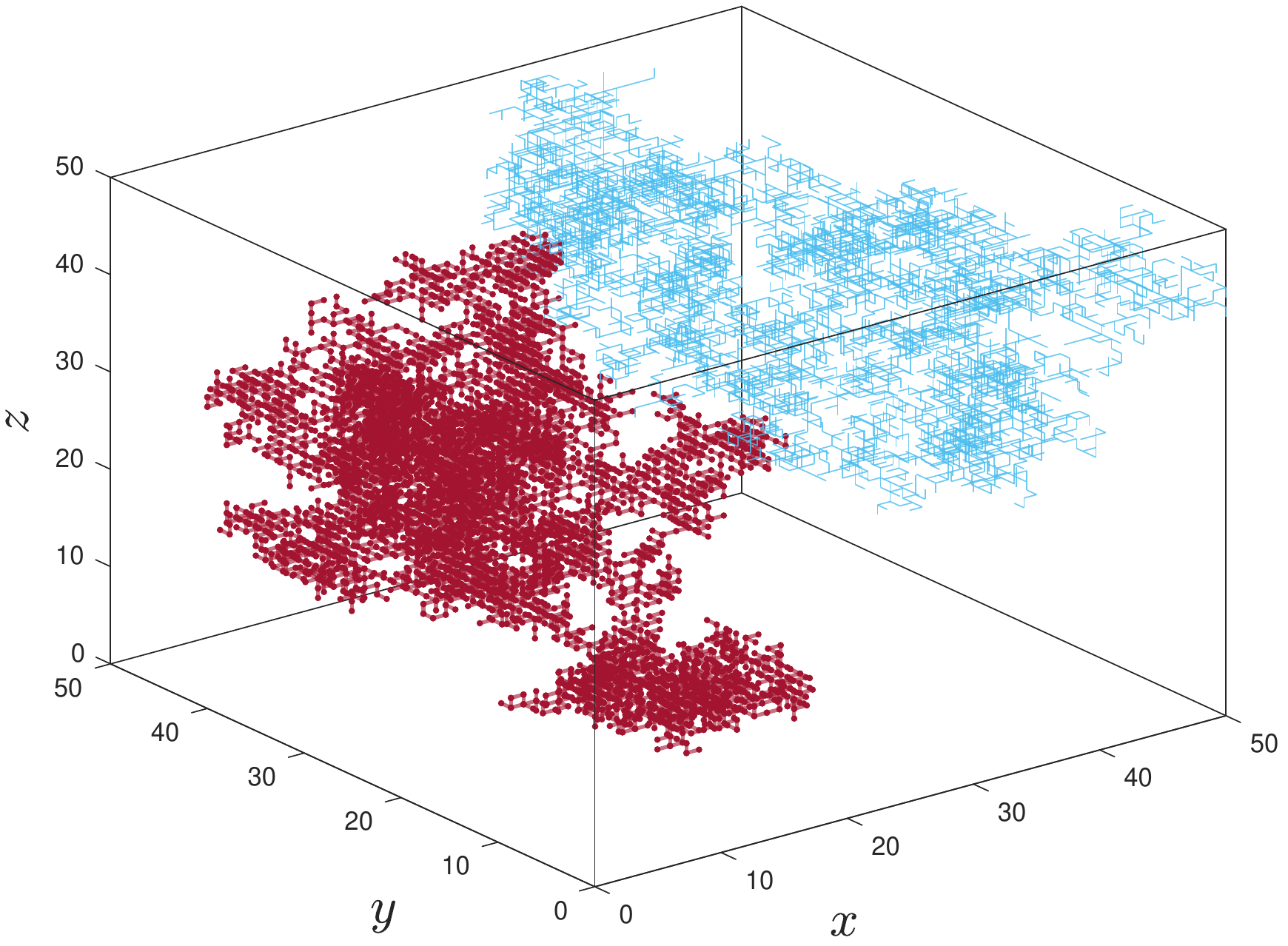}\label{fig_1b}}\\
\centering
\subfigure[T = 2.0]{\includegraphics[width=0.49\textwidth,height=0.30\textwidth]{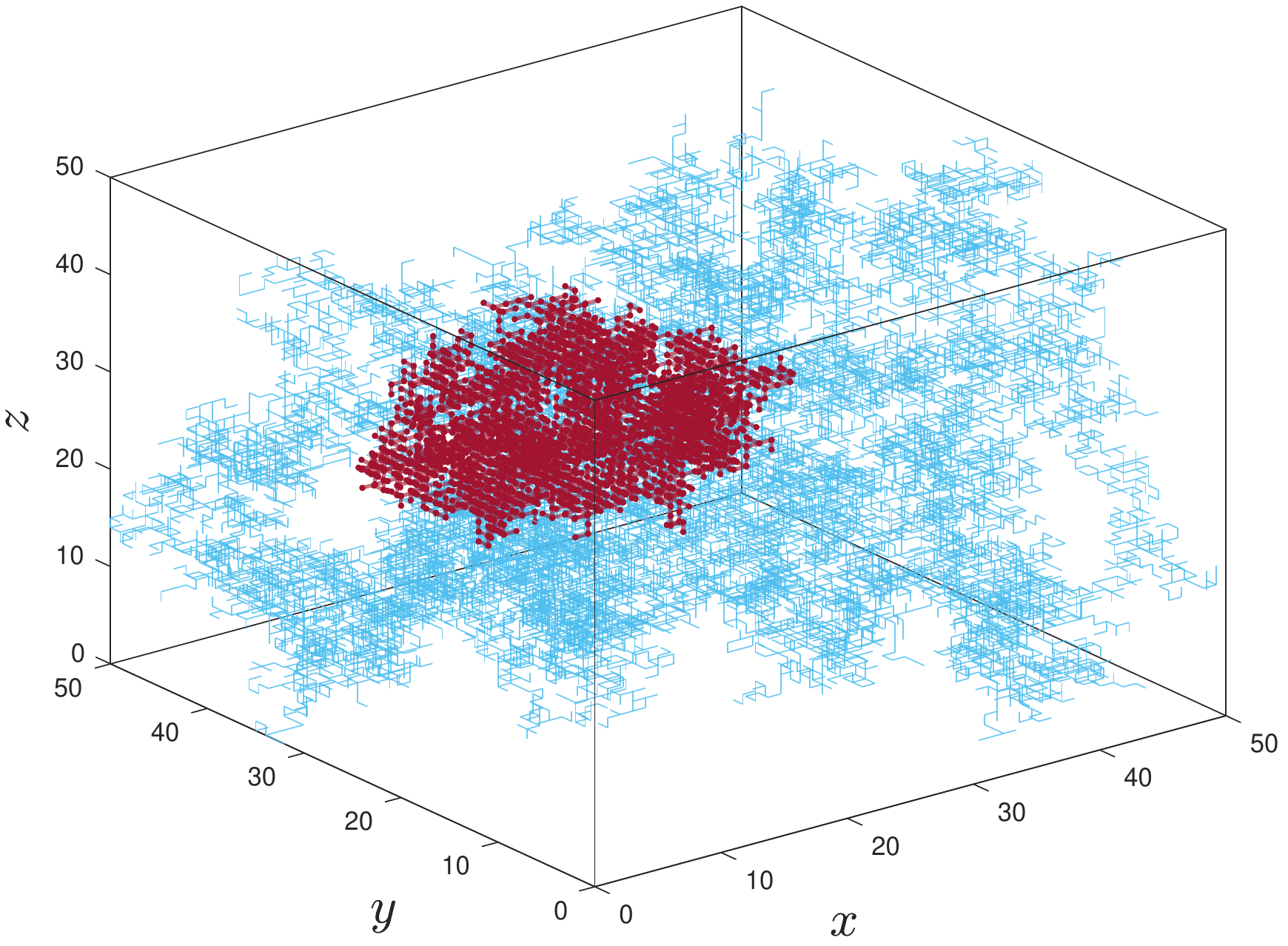}\label{fig_1c}}
\subfigure[T = 2.2]{\includegraphics[width=0.49\textwidth,height=0.30\textwidth]{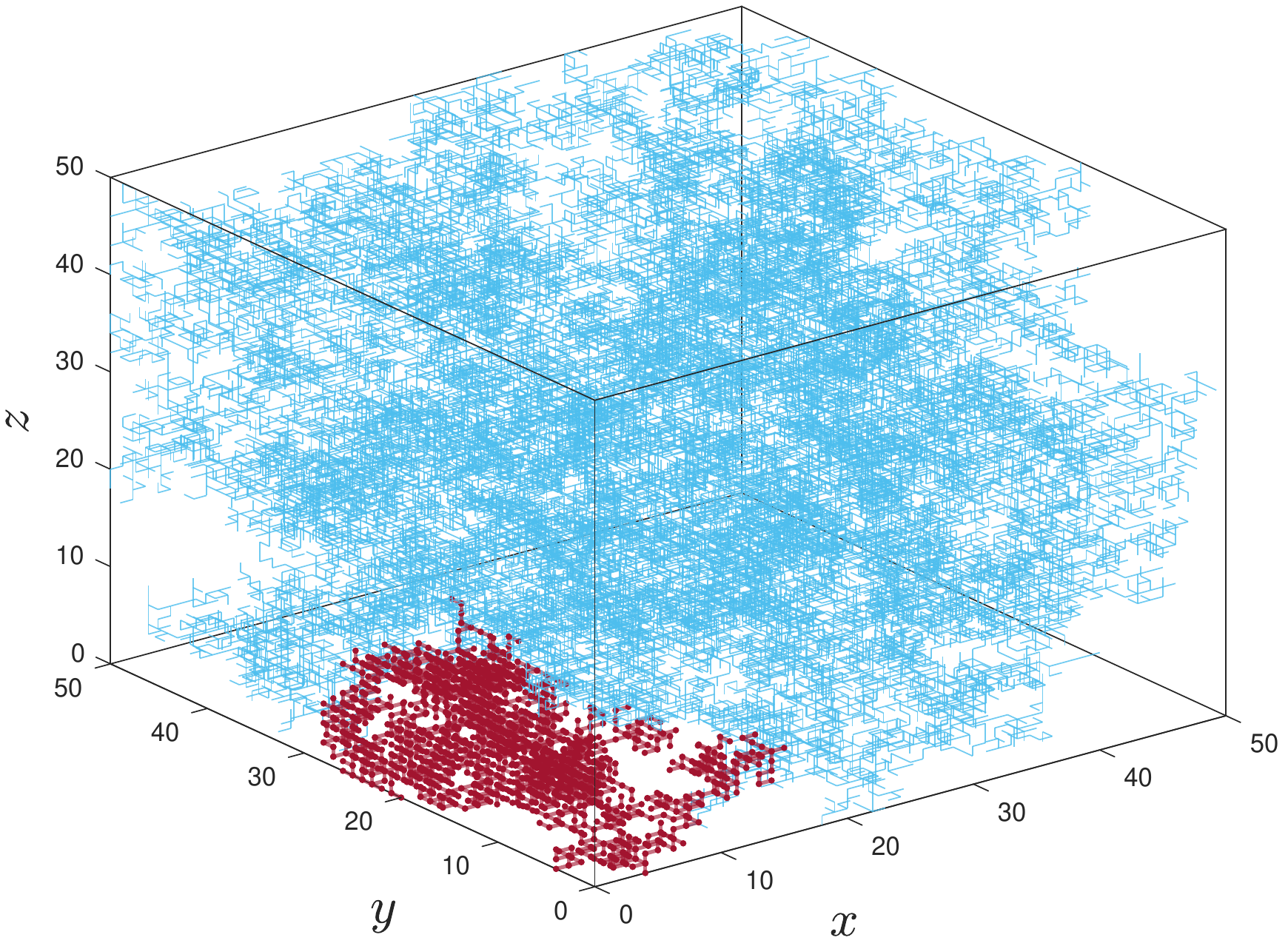}\label{fig_1d}}
%\vspace{2.5 in}
\caption{Illustrations of FRI truncated in a box $[0, 50]^3$: The first and second largest clusters under $u = 1/6$ and different parameter {$T's$}. The sub and super critical phases are demonstrated under $T= 1.4$ and $T=2.2$, respectively. The intermediate phases under $T=1.8$ and $T=2.0$ provide some evidence on the percolation phase transition.}
\label{fig_example}
\end{figure}

In the followings we denote by $\mathcal{FI}_d^{u,T}$ the FRI in $\BZ^d$ with multiplicative parameter $u$ and truncation parameter $T$. A key character of the FRI is percolation property, i.e., the existence and uniqueness of an infinite cluster within $\mathcal{FI}_d^{u,T}$. In contrast to RI, where $\mathcal{I}_d^u$ almost surely percolates for all $d$ and $u>0$, FRI has been proved in \cite{FRI_1} to have the following phase transition as an edge percolation model:
\begin{itemize}
\item Supercritical phase (Theorem 1,  \cite{FRI_1}): 	for all $d\ge 3$ and $u >0$, there is a $0 < T_1 (u,d) < \infty$ such that for all $T > T_1$, $\mathcal{FI}_d^{u,T}$ has a unique infinite cluster almost surely.
\item Subcritical phase (Theorem 2, \cite{FRI_1}): for all $d\ge 3$ and $u >0$, there is a $0< T_0 (u,d) < \infty$ such that for all $0 < T < T_0$, $\mathcal{FI}_d^{u,T}$ has no infinite cluster almost surely. 
\end{itemize} 
Intuitively, the percolation can be visualized by running one realization under different parameters $u, T$ and plotting the first and second largest clusters restricted in a finite box $[0, 50]^3$. Two small clusters in Fig.~\ref{fig_1a}, that corresponds to $u=1/6, T=1.4$,  provide some evidence that no infinite cluster exists, while a huge cluster along with a smaller cluster in Fig.~\ref{fig_1d}, that corresponds to $u=1/6, T=2.2$, indicates that there may exist only one infinite cluster.  Intermediate phases between the subcritical one and supercritical one are presented in Figs.~\ref{fig_1b} and \ref{fig_1c}. One can see that the phase transition may occur near $T = 1.8$, in which the first and second largest clusters are almost of the same size. The dominance of the first largest cluster can be apparently observed as $T$ goes larger, say, $T=2.0$.

%\note{3 subfigures of the largest+ 2nd}

Moreover, a follow-up work \cite{FRI_2} proved recently that for all $d\ge 3$ and $u >0$, there is a $T_2 (u,d) \in [T_1(u,d) , \infty)$ such that for all $T > T_2$, the chemical distance on $\mathcal{FI}_d^{u,T}$ is asymptotically of the same order as the Euclidean distance.  \cite{FRI_2} further proves that FRI has local uniqueness property for all sufficiently large $T$. See \cite{drewitz2014chemical}, for precise definitions for chemical distance and local uniqueness. However, since $\mathcal{FI}_d^{u,T}$ may be non-monotonic with respect to $T$, the existence of a subcritical and a supercritical phase is insufficient to guarantee the sharpness of the phase transition. It is conjectured in  \cite{FRI_1} that there is a unique critical value $T_c(u,d)$ such that $\mathcal{FI}_d^{u,T}$ percolates when $T>T_c$ and has no infinite cluster almost surely when $T<T_c$.

The percolation phase transition is closely related to the tradeoff mechanism with respect to the parameter $T$: As $T$ increases, there will be in average fewer and fewer fibers starting from each vertex. But in compensation, the length of each remaining fiber increases proportional to $T$, so that we are less likely to see the start or end of any fiber locally. In fact as $T\to \infty$, an FRI increasingly resembles the limiting model (which is the classical RI itself) where all fibers are doubly infinite SRW trajectories. And it has been proved in \cite{bowen2019finitary} that $\mathcal{FI}_d^{u,T}\Rightarrow \mathcal{I}_d^u$ under the weak-* topology. With the observation above, it is natural to ask how the FRI evolves with respect to $T$, or more specifically, whether or not it has stochastic monotonicity \cite{FRI_1}.

In this paper, we prove that for all $u>0$, $\mathcal{FI}_d^{u,T}$ is NOT monotonically increasing in $T$ for $d=3,4$. This finding reveals the evolution of FRI with respect to $T$ might be more nontrivial than previous thought, and makes the characterization of the lower dimensional phase transition in $T$ a more interesting question. At the same time, as $d$ goes to infinity, our calculation leads us to doubt that stochastic monotonicity may become true when $d$ is sufficiently high. 

The non-monotonicity of low dimensional FRI casts shadows on the conjecture of existence and uniqueness of critical value $T_c$. A large-scale parallel computing algorithm is employed to explore the behavior of FRI when the fiber length factor $T$ changes in the intermediate phase $[T_0,T_1]$. Our numerical findings strongly suggest that, although no longer monotonic, for any $d\ge 3, \ u>0$, $\mathcal{FI}_d^{u,T}$ still has a sharp and unique phase transition of percolation. Our numerical tests also indicate that the phase diagram is estimated to be an inversely proportional function with respect to $u$. We also rigorously prove that there is some constant $C<\infty$ such that for all sufficiently large $u$, $T_1\le C u^{-1}$, which also suggests that it is unlikely to have multiple phase transitions, especially when $u$ is large. 

The rest of this paper is outlined as follows: in Section \ref{sec_def}, we recall the precise definition of FRI, together with some important notations and constructions crucial for our theoretical and numerical discussions; in Section \ref{sec_theo} we discuss the non-monotonicity of FRI for $d=3, 4$; while our numerical explorations on the phase diagram are presented in Section \ref{sec_num}.

%%%%%%%%%%%%%%%%%%%%%%%%%%%%%%%%%%%%%%%%%%
%\section{Results}

\section{Definitions and Notations}
\label{sec_def}

In this section, we recall the precise definition of FRI, together with some important notations and constructions in \cite{bowen2019finitary, FRI_1, FRI_2}. We start with some standard notations for simple random walks. Without causing further confusion, we will use $\BZ^d$ to denote both vertices and (the nearest neighbor) edges in the $d-$dimensional lattice throughout this paper. Then for a subgraph $\mathcal{G}=(V,E)\subset \BZ^d$, we call it connected if any $v_1,v_2\in V$ can be connected by a collection of edges in $E$. And for $1\le j\le d$, let $x_j\in \mathbb{R}^d$ satisfy $x_j^{(i)}=\mathbbm{1}_{\{i=j\}}$, $1\le i\le d$. Denote the edges $\{0,x_j\}$ and $\{0,-x_j\}$ by $e_j$ and $-e_j$. Note that $\{e_j\}_{j=1}^d$ form the basis of $\BZ^d$. For any subset of vertices $A\subset \BZ^d$, let
\begin{align*}
&\partial^{in}A=\left\{x\in A, \ s.t. \exists y\in A^c, \ \|x-y\|_1=1 \right\}\\
&\partial^{out}A=\left\{x\in A^c, \ s.t. \exists y\in A, \ \|x-y\|_1=1 \right\}
\end{align*}
be the inner and outer boundary of $A$. 

In this paper, we denote by $\{X_n\}_{n=0}^\infty$ a simple random walk (SRW) in $\BZ^d$ starting from $X_0$, with its distribution denoted by $P_{X_0}(\cdot)$. Note that for any integer $0\le n_0\le \infty$, the SRW trajectory $\{X_n\}_{n=0}^{n_0}$ naturally induces a collection of edges $\big\{ \{X_i, X_{i+1}\}\big\}_{i=0}^{n_0-1}$. Moreover, for any $T>0$, let $G_T\ge 0$ be a geometric random variable with $p=1/(T+1)$ which is independent to $\{X_n\}_{n=0}^\infty$. Then we call $\{X_n\}_{n=0}^{G_T}$ a geometrically killed SRW with parameter $T$, and denote its distribution by $P_{X_0}^{(T)}(\cdot)$ with the convention $P_{X_0}^{(\infty)}(\cdot)=P_{X_0}(\cdot)$. 

Moreover, we denote by 
\begin{align*}
&\bar H_{d,A}^{(T)}=\inf\{n\ge 0, \ X^{(T)}_n\in A\}\\
& H_{d,A}^{(T)}=\inf\{n\ge 1, \ X^{(T)}_n\in A\}
\end{align*}
the first hitting and first returning times to $A$, with the convention $\inf \emptyset =\infty$. 
\begin{Remark}
It is worth to note that  $H_{d,A}^{(T)}\equiv \infty$ when $G_T=0$. 
\end{Remark}

For a finite subset $A$ and vertex $x\in \BZ^d$, define  {\bf killed escape probability}
\begin{equation}
\label{eq_escape}
Es^{(T)}_{d,A}(x)=P_x^{(T)}(H_{d,A}^{(T)}=\infty),
\end{equation} 
together with equilibrium measure
\begin{equation}
\label{eq_equilibrium}
e^{(T)}_{d,A}(x)=(2d)\times Es^{(T)}_{d,A}(x) \ind_{x\in A}
\end{equation}
and killed capacity
\begin{equation}
\label{eq_capacity}
\text{cap}_d^{(T)}(A)=\sum_{x\in \BZ^d} e^{(T)}_{d,A}(x). 
\end{equation}

We also abbreviate $\bar H_{d,A}^{(T)}$, $H_{d,A}^{(T)}$, $Es^{(T)}_{d,A}$, $e^{(T)}_{d,A}$ and $\text{cap}_d^{(T)}$ to $\bar H_{d,A}$, $H_{d,A}$, $Es_{d,A}$, $e_{d,A}$ and $\text{cap}_d$ respectively, when $T=\infty$.

\begin{Remark}
The factor $2d$ in killed equilibrium measure in added for the technical reason to make FRI converges to RI with the same intensity as $T\to\infty$. See \cite{bowen2019finitary} for details. 
\end{Remark}

It is worth to note that for $x\in A\setminus \partial^{in} A$, $e^{(T)}_{d,A}(x)=P(G_T=0)=1/(1+T)$. 

\subsection{Definition of FRI}

According to \cite{FRI_1}, there are two equivalent definitions of FRI. Denote the set of all finite paths on $\mathbb{Z}^d$ by $W_d^{\left[ 0,\infty\right) }$. Since $W_d^{\left[ 0,\infty\right) }$ is countable, the measure $v_d^{(T)}=\sum_{x\in \mathbb{Z}^d}\frac{2d}{T+1}P_x^{(T)}$ is a $\sigma-$ finite measure on $W_d^{\left[ 0,\infty\right) }$.
\begin{Definition}\label{definition1}
	For $0<u,T<\infty$, the finitary random interlacements $\mathcal{FI}^{u,T}_d$ is a Poisson point process on $W_d^{\left[ 0,\infty\right) }$ with intensity measure $uv_d^{(T)}$. The law of $\mathcal{FI}_d^{u,T}$ is denoted by $P^{u,T}$. 
\end{Definition}
\begin{Definition}\label{definition2}
	For each site $x\in \mathbb{Z}^d$, $N_x$ is a Poisson random variable with parameter $\frac{2du}{T+1}$. Start $N_x$ independent geometrically killed simple random walks starting at $x$ with killing rate $\frac{1}{T+1}$. Then one may equivalently define  $\mathcal{FI}_d^{u,T}$ as a point measure on $W_d^{\left[ 0,\infty\right) }$ composed of all the trajectories above from all sites in $\mathbb{Z}^d$. 
\end{Definition}

\subsection{Configurations within a finite set}\label{2.2}
In fact, given a finite set $K\subset \mathbb{Z}^d$, the distribution of FRI within $K$ can be discribed precisely. By Lemma 2.2 of \cite{FRI_1}, if we start $N_x\sim Pois\left( u*e_{d,K}^{(T)}(x)\right)$ independent random walks with distribution $P_x^{(T)}$ for any $x\in K$ (denote all these trajectories by $\left\lbrace\eta_i\right\rbrace_{i=1}^{N_K} $), then $\bigcup_{i=1}^{N_K}\eta_i\cap K$ has the same distribution as $\mathcal{FI}_d^{u,T}\cap K$.

\subsection{Stochastic dominance and monotonicity}
A sufficient condition for the existence of the critical value $T_c$ is the stochastic monotonicity with respect to $T$. More precisely, if for any $T'>T$, there is a coupling between $\mathcal{FI}_d^{u,T'}$ and $\mathcal{FI}_d^{u,T}$ such that $\mathcal{FI}_d^{u,T}\subset\mathcal{FI}_d^{u,T'}$ almost surely, then $T_c$ must exist. Therefore, we need the concept of stochastic dominance to describe the existence of the coupling. 
\begin{Definition}(Definition 2.1, Chapter 2 of \cite{liggett2012interacting})
	Assume that $X$ is a compact metric space with a given partial order. Say a function $f$ on $X$ is monotone if $f(\eta)\le f(\zeta)$ for any $\eta,\zeta\in X$, $\eta\le \zeta$. Then for two probability measures $\mu_1,\mu_2$ on $X$, say $\mu_2$ stochastically dominates $\mu_1$ (written by $\mu_1\overset{d}{\le}\mu_2$) if and only if for any monotone function $f$ on $X$, $$\int fd\mu_1\le \int fd\mu_2.$$
\end{Definition}
By Theorem 2.4 in the Chapter 2 of \cite{liggett2012interacting}, we know that the coupling mentioned above exists if and only if $\mathcal{FI}_d^{u,T}\overset{d}{\le}\mathcal{FI}_d^{u,T'}$.
\section{Single edge density and monotonicity}
\label{sec_theo}

Recalling the definition  of stochastic monotonicity, to prove/disprove $\mathcal{FI}_d^{u,T_1} \overset{d}{\ge} \mathcal{FI}_d^{u,T_2}$ for all $T_1\ge T_2$, it is equivalent to verify whether or not for all monotonically increasing function $\varphi$ on $\BZ^d$, one always has
$$
E\left[\varphi\left(\mathcal{FI}_d^{u,T_1} \right) \right]\ge E\left[\varphi\left(\mathcal{FI}_d^{u,T_2} \right) \right]. 
$$
In particular, we can first take the test function as the very simple form $\varphi_1(E)=\ind_{e_1\in E}$ and thus 
\begin{equation}
E\left[\varphi_1\left(\mathcal{FI}_d^{u,T} \right) \right]=P\left(e_1\in \mathcal{FI}_d^{u,T}\right)\overset{\Delta}{=}p_{d,u}(T)
\end{equation}
gives the probability that any single (undirected) edge is traversed by the FRI. 

\begin{Theorem}
\label{d34}
For $d=3,4$ and any $u\in (0,\infty)$, $p_{d,u}(T)\in C^1(0,\infty)$. Moreover, there are $t_0>0$ and $T_0<\infty$ such that 
\begin{itemize}
\item $p_{d,u}'(T)>0$ for all $T\in (0,t_0)$;
\item $p_{d,u}'(T)<0$ for all $T\in (T_0,\infty)$.
\end{itemize}
Thus $\mathcal{FI}_d^{u,T} $ is NOT stochastically monotone with respect to $T$ for $d=3,4$ and any $u>0$. 
\end{Theorem}

Before presenting the proof of Theorem \ref{d34}, we first cite the following useful result from \cite{lawler2010random} on the expected length of excursion given a SRW returns to where it starts. 

\begin{Proposition}[Exercise 11.1, \cite{lawler2010random}]
\label{Lawler_11_1}
Suppose $d\ge 3$ and $Y_n$ is a simple random walk in $\BZ^d$ conditioned on return to the origin. Then 
%\begin{enumerate}[(i)]
\begin{enumerate}
\item For all $d\ge 3$, $Y_n$ is a recurrent Markov chain. 
\item Assume $Y_0=0$ and let $\Gamma=\min\{j>0: \ Y_j=0\}$. Then
$$
P(\Gamma=2n)\asymp n^{-d/2}, \ n\to\infty.
$$ 
In particular,
$$
E[T]\left\{
\begin{aligned}
=\infty, \ d\le 4,\\
<\infty, \ d\ge 5.
\end{aligned}
\right.
$$
\end{enumerate}
\end{Proposition}

\begin{proof}[Proof of Theorem \ref{d34}]
	
	We denote the event $\{e_1\ is\ not\ contained\ by\ the\ trajectory\}$ by $F$ and then we are going to calculate $P^{u,T}\left(e_1\notin \mathcal{FI}_d^{u,T}\right)$. First, we need to calculate $P_{x_2}^{(T)}\left(F\right)$ and $P_{-x_1}^{(T)}\left(F\right)$. Denote that $E_1=P_{-x_1}^{(T)}\left(H^{(T)}_{d,\{0,x_1\}}=\infty\right) $ and $E_2=P_{x_2}^{(T)}\left(H^{(T)}_{d,\{0,x_1\}}=\infty\right) $. We have 
	\begin{equation}\label{2}
		\begin{split}
			P_{x_2}^{(T)}\left(F\right)=&E_2+\sum\limits_{n=1}^{\infty}P_{x_2}^{(T)}\left(H^{(T)}_{d,\{0,x_1\}}=n, F \right) \\
			=&E_2+\sum\limits_{n=1}^{\infty}P_{x_2}^{(T)}\left(H^{(T)}_{d,\{0,x_1\}}=n \right)\left[\frac{1}{T+1}+\frac{T}{T+1}\left( \frac{2d-2}{2d}P_{x_2}^{(T)}\left(F\right)+\frac{1}{2d}P_{-x_1}^{(T)}\left(F\right)\right) \right]\\
			=&E_2+\left(1-E_2 \right) \left[\frac{1}{T+1}+\frac{T}{T+1}\left( \frac{2d-2}{2d}P_{x_2}^{(T)}\left(F\right)+\frac{1}{2d}P_{-x_1}^{(T)}\left(F\right)\right) \right].
		\end{split}
	\end{equation}
	%Thus, \begin{equation}
	%	P_{x_2}^{(T)}\left(F\right)\left[1-\frac{T}{T+1}*\frac{2d-2}{2d}*\left(1-E_2  \right)  \right] =\left(1-E_2 \right)*\frac{T}{T+1}*\frac{1}{2d}P_{-x_1}^{(T)}\left(F\right)+E_2+\left(1-E_2 \right)*\frac{1}{T+1}.
	%\end{equation}
	In the same way, we have \begin{equation}\label{3}
		P_{-x_1}^{(T)}\left(F\right)=E_1+\left(1-E_1 \right) \left[\frac{1}{T+1}+\frac{T}{T+1}\left( \frac{2d-2}{2d}P_{x_2}^{(T)}\left(F\right)+\frac{1}{2d}P_{-x_1}^{(T)}\left(F\right)\right) \right].
	\end{equation}
	%Thus, \begin{equation}\label{4}
	%	P_{-x_1}^{(T)}\left(F\right)\left[1-\left(1-E_1 \right)*\frac{T}{T+1}*\frac{1}{2d} \right]=\left(1-E_1 \right)*\frac{T}{T+1}*\frac{2d-2}{2d}P_{x_2}^{(T)}\left(F\right)+ E_1+\left(1-E_1 \right)\frac{1}{T+1}.
	%\end{equation}
	Combine (\ref{2}) and (\ref{3}), 
	\begin{equation}\label{5}
		P_{x_2}^{(T)}\left(F\right)\left[1+(E_2-E_1)*\frac{T}{T+1}*\frac{2d-2}{2d}\right] 	=P_{-x_1}^{(T)}\left(F\right)\left[1+(E_1-E_2)\frac{T}{T+1}*\frac{1}{2d} \right] +(E_2-E_1)\frac{T}{T+1}.
	\end{equation}
	By (\ref{2}) and (\ref{5}), we have \begin{equation}
		\begin{split}
			&P_{-x_1}^{(T)}\left(F\right)*\frac{\left(1-E_2 \right)*\frac{T}{T+1}*\frac{1}{2d}}{1-\frac{T}{T+1}*\frac{2d-2}{2d}*\left(1-E_2  \right)}+\frac{E_2+\left(1-E_2 \right)*\frac{1}{T+1}}{1-\frac{T}{T+1}*\frac{2d-2}{2d}*\left(1-E_2  \right)}\\
			=&P_{-x_1}^{(T)}\left(F\right)*\frac{1-(E_2-E_1)*\frac{T}{T+1}*\frac{1}{2d}}{1+(E_2-E_1)*\frac{T}{T+1}*\frac{2d-2}{2d}}+\frac{(E_2-E_1)\frac{T}{T+1}}{1+(E_2-E_1)*\frac{T}{T+1}*\frac{2d-2}{2d}}.
		\end{split}
	\end{equation}
	Therefore, \begin{equation}\label{7}
		\begin{split}	
			%		=&\frac{\left[1-\frac{T}{T+1}*\frac{2d-2}{2d}*\left(1-E_2  \right) \right]*\left[(E_2-E_1)\frac{T}{T+1} \right]-\left[E_2+\left(1-E_2 \right)*\frac{1}{T+1}\right]*\left[1+(E_2-E_1)*\frac{T}{T+1}*\frac{2d-2}{2d} \right]    }{\left[\left(1-E_2 \right)*\frac{T}{T+1}*\frac{1}{2d} \right]*\left[1+(E_2-E_1)*\frac{T}{T+1}*\frac{2d-2}{2d} \right] -\left[1-\frac{T}{T+1}*\frac{2d-2}{2d}*\left(1-E_2  \right) \right]*\left[1-(E_2-E_1)*\frac{T}{T+1}*\frac{1}{2d} \right]  }\\
			P_{-x_1}^{(T)}\left(F\right)=&\frac{E_2+(1-E_2)\frac{1}{T+1}-(E_2-E_1)\frac{2}{2d}\frac{T}{T+1}}{1-\left[\frac{2d-1}{2d}(1-E_2)+\frac{1}{2d}(E_2-E_1) \right]\frac{T}{T+1}}\\
			=&\frac{\left[\frac{2d-2}{2d}E_2+\frac{2}{2d}E_1 \right]T+1 }{\left[\frac{2d-2}{2d}E_2+\frac{1}{2d}E_1+\frac{1}{2d} \right]T+1 }.
		\end{split}
	\end{equation}
	In the same way, we have \begin{equation}
		\begin{split}
			P_{x_2}^{(T)}\left(F\right)=\frac{\left[\frac{2d-1}{2d}E_2+\frac{1}{2d}E_1 \right]T+1 }{\left[\frac{2d-2}{2d}E_2+\frac{1}{2d}E_1+\frac{1}{2d} \right]T+1 }.
			%	&\frac{E_2+(1-E_2)*\frac{1}{T+1}+(E_1-E_2)*\frac{1}{2d}*\frac{T}{T+1}}{1-\frac{T}{T+1}\frac{2d-2}{2d}(1-E_2)-\frac{T}{T+1}\frac{1}{2d}(1-E_1)}\\	\frac{\frac{E_2+(1-E_2)*\frac{1}{T+1}}{(1-E_2)*\frac{T}{T+1}*\frac{1}{2d}}-\frac{(E_2-E_1)*\frac{T}{T+1}}{1+(E_1-E_2)*\frac{T}{T+1}*\frac{1}{2d}}}{\frac{1-\frac{T}{T+1}*\frac{2d-2}{2d}*(1-E_2)}{(1-E_2)*\frac{T}{T+1}*\frac{1}{2d}}-\frac{1+(E_2-E_1)*\frac{T}{T+1}*\frac{2d-2}{2d}}{1+(E_1-E_2)*\frac{T}{T+1}*\frac{1}{2d}}}\\
			%		=&\frac{\left[E_2+(1-E_2)*\frac{1}{T+1} \right]*\left[1+(E_1-E_2)*\frac{T}{T+1}*\frac{1}{2d} \right]-\left[(E_2-E_1)*\frac{T}{T+1} \right]*\left[(1-E_2)*\frac{T}{T+1}*\frac{1}{2d} \right]}{\left[1-\frac{T}{T+1}*\frac{2d-2}{2d}*(1-E_2) \right]*\left[1+(E_1-E_2)*\frac{T}{T+1}*\frac{1}{2d} \right]-\left[1+(E_2-E_1)*\frac{T}{T+1}*\frac{2d-2}{2d} \right]*\left[(1-E_2)*\frac{T}{T+1}*\frac{1}{2d} \right] }\\
		\end{split}
	\end{equation}
	Therefore,  \begin{equation}
		\begin{split}
			P_{0}^{(T)}\left(F\right)=&\frac{1}{T+1}+\frac{T}{T+1}\left[\frac{2d-2}{2d}P_{x_2}^{(T)}\left(F\right)+\frac{1}{2d}P_{-x_1}^{(T)}\left(F\right) \right]\\
			%		=&\frac{1}{T+1}+\frac{T}{T+1}*\frac{\left[\frac{2d-2}{2d}E_2+\frac{1}{2d}E_1 \right]T+\frac{2d-1}{2d} }{\left[\frac{2d-2}{2d}E_2+\frac{1}{2d}E_1+\frac{1}{2d} \right]T+1 }\\
			%		=&\frac{\left[\frac{2d-2}{2d}E_2+\frac{1}{2d}E_1 \right]T+\frac{2d-1}{2d} }{\left[\frac{2d-2}{2d}E_2+\frac{1}{2d}E_1+\frac{1}{2d} \right]T+1 }+\frac{1}{T+1}\frac{\frac{1}{2d}(T+1) }{\left[\frac{2d-2}{2d}E_2+\frac{1}{2d}E_1+\frac{1}{2d} \right]T+1 }\\
			=&\frac{\left[\frac{2d-2}{2d}E_2+\frac{1}{2d}E_1 \right]T+1 }{\left[\frac{2d-2}{2d}E_2+\frac{1}{2d}E_1+\frac{1}{2d} \right]T+1 }.
		\end{split}
	\end{equation}
	Restricted on $\{0,x_1\}$, there are $ Pois\left(2du*Es_{d,\{0,x_1\}}^{(T)}(0)\right) $ independent trajectories starting from $0$ and $ Pois\Big(2du*Es_{d,\{0,x_1\}}^{(T)}(1)\Big) $ trajectories starting from $1$. Note that $Es_{d,\{0,x_1\}}^{(T)}(0)=Es_{d,\{0,x_1\}}^{(T)}(1)$, we have 
	\begin{equation}
		\begin{split}
			P^{u,T}\left(e_1\notin \mathcal{FI}_d^{u,T} \right) =&\left[ \sum_{m=0}^{\infty}\exp(-2du*Es_{d,\{0,x_1\}}^{(T)}(0))\frac{\left( 2du*Es_{d,\{0,x_1\}}^{(T)}(0)\right)^m}{m!}\left( P_{0}^{(T)}\left(F\right)\right)^m \right]^2\\
			=&\exp(-4du*Es_{d,\{0,x_1\}}^{(T)}(0)\left(1-P_{0}^{(T)}\left(F\right) \right) ).
		\end{split}
	\end{equation}
	%Therefore, $p_{d,u}(T)=1-P^{u,T}\left(e_1\notin \mathcal{FI}_d^{u,T} \right)\in C^1(0,\infty) $.

	%it's sufficient to check that whether $Es_{\{0,x_1\}}^{(T)}(0)\left(1-P_{0}^{(T)}\left(F\right) \right) $ has monotonicity. We will calculate $(f*g)'(T)$.
	
	Let $f(T)=Es_{d,\{0,x_1\}}^{(T)}(0)=Es_{d,\{0,x_1\}}(0)+\sum\limits_{n=1}^{\infty}P_0\left(H_{d,\{0,x_1\}}=n \right)\left(1-\left(1-\frac{1}{T+1} \right)^n  \right)$ and $g(T)=1-P_{0}^{(T)}\left(F\right)=\frac{\frac{1}{2d}T }{\left[\frac{2d-2}{2d}E_2+\frac{1}{2d}E_1+\frac{1}{2d} \right]T+1 }$. We have \begin{equation}\label{310}
		\begin{split}
			f'(T)=&-\frac{1}{(T+1)^2}\sum\limits_{n=1}^{\infty}\left[P_0\left(H_{d,\{0,x_1\}}=n \right)*n*\left(1-\frac{1}{T+1} \right)^{n-1} \right]\\
			=&-\frac{1}{T(T+1)}\sum\limits_{n=1}^{\infty}\left[P_0\left(H_{d,\{0,x_1\}}=n \right)*n*\left(1-\frac{1}{T+1} \right)^{n}\right]\\
			=&-\frac{1}{T(T+1)}*E_0^{(T)}\left[H^{(T)}_{d,\{0,x_1\}}; 1\le H^{(T)}_{d,\{0,x_1\}}<\infty \right]<\infty.
		\end{split}
	\end{equation}
	Meanwhile, \begin{equation}\label{311}
		\begin{split}
			g'(T)=\frac{\frac{1}{2d} }{\left[ \left(\frac{2d-2}{2d}E_2+\frac{1}{2d}E_1+\frac{1}{2d} \right)+\frac{1}{T} \right]^2}*\frac{1}{T^2}=\frac{\frac{1}{2d} }{\left[ \left(\frac{2d-2}{2d}E_2+\frac{1}{2d}E_1+\frac{1}{2d} \right)T+1 \right]^2}\overset{\Delta}{=}\frac{\frac{1}{2d}}{\left(aT+1 \right)^2 },
		\end{split}
	\end{equation}
	where $a=\frac{2d-2}{2d}E_2+\frac{1}{2d}E_1+\frac{1}{2d}$. Combine (\ref{310}) and (\ref{311}),
	\begin{equation}\label{312}
		\begin{split}
			(f\cdot g)'(T)=&Es^{(T)}_{d,\{0,x_1\}}(0)*\frac{\frac{1}{2d}}{\left(aT+1 \right)^2 }-\frac{1}{T(T+1)}*E_0^{(T)}\left[ H^{(T)}_{d,\{0,x_1\}}; 1\le H^{(T)}_{d,\{0,x_1\}}<\infty \right]*\frac{\frac{1}{2d}T}{aT+1}\\
			=&\frac{\frac{1}{2d}}{(aT+1)^2}\left[Es^{(T)}_{d,\{0,x_1\}}(0)-\frac{aT+1 }{T+1}*E_0^{(T)}\left[  H^{(T)}_{d,\{0,x_1\}}; 1\le H^{(T)}_{d,\{0,x_1\}}<\infty \right]  \right]<\infty.
		\end{split}
	\end{equation}
	Therefore, $p_{d,u}(T)=1-\exp(-4du*f(T)*g(T))\in C^1(0,\infty)$.
	
	Note that $\lim\limits_{T\to 0+}Es^{(T)}_{d,\{0,x_1\}}(0)=1$ and $\lim\limits_{T\to 0+}E_0^{(T)}\left[H^{(T)}_{d,\{0,x_1\}}; 1\le H^{(T)}_{d,\{0,x_1\}}<\infty \right]=0 $, we have\begin{equation}
		\lim\limits_{T\to 0+}(f\cdot g)'(T)=\frac{1}{2d}>0.
	\end{equation} 
	Similar to Proposition \ref{Lawler_11_1}, for $d=3,4$, $E_0\left[  H_{\{0,x_1\}}; 1\le H_{\{0,x_1\}}<\infty \right]=\infty$. Then we have \begin{equation}
		\lim\limits_{T\to \infty}E_0^{(T)}\left[  H^{(T)}_{d,\{0,x_1\}}; 1\le H^{(T)}_{d,\{0,x_1\}}<\infty \right]=\infty.
	\end{equation}
	Note that $\lim\limits_{T\to \infty}Es^{(T)}_{d,\{0,x_1\}}(0)=Es_{d,\{0,x_1\}}(0)$, we know that $\exists T_0>0$ such that $\forall T>T_0$,\begin{equation}
		Es^{(T)}_{d,\{0,x_1\}}(0)-\frac{aT+1 }{T+1}*E_0^{(T)}\left[  H^{(T)}_{d,\{0,x_1\}}; 1\le H^{(T)}_{d,\{0,x_1\}}<\infty \right] <0.
	\end{equation}
	
	In conclusion, we know that when $T$ is small, $(f\cdot g)'(T)>0$ and when $T$ is large, $(f\cdot g)'(T)<0$. Recalling that $p_{d,u}(T)=1-\exp(-4du*f(T)*g(T))$, the proof is complete. 
	% Therefore, $Es_{\{0,x_1\}}^{(T)}(0)\left(1-P_{0}^{(T)}\left(F\right) \right) $ doesn't have monotonicity and so does $P^{u,T}\left(e_1\notin \mathcal{FI}_d^{u,T} \right)$.
\end{proof}

\begin{Remark}
	Though it is true that Proposition \ref{Lawler_11_1} as stated is for the expected time of returning to 0 rather than to $\{0,x_1\}$, the result and proof are exactly parallel for returning to any finite subset. So we decide to cite \cite{lawler2010random} rather than repeat the proof. 
\end{Remark}

\begin{figure}[H]
\centering
\includegraphics[width=4.8in,height=3.2in]{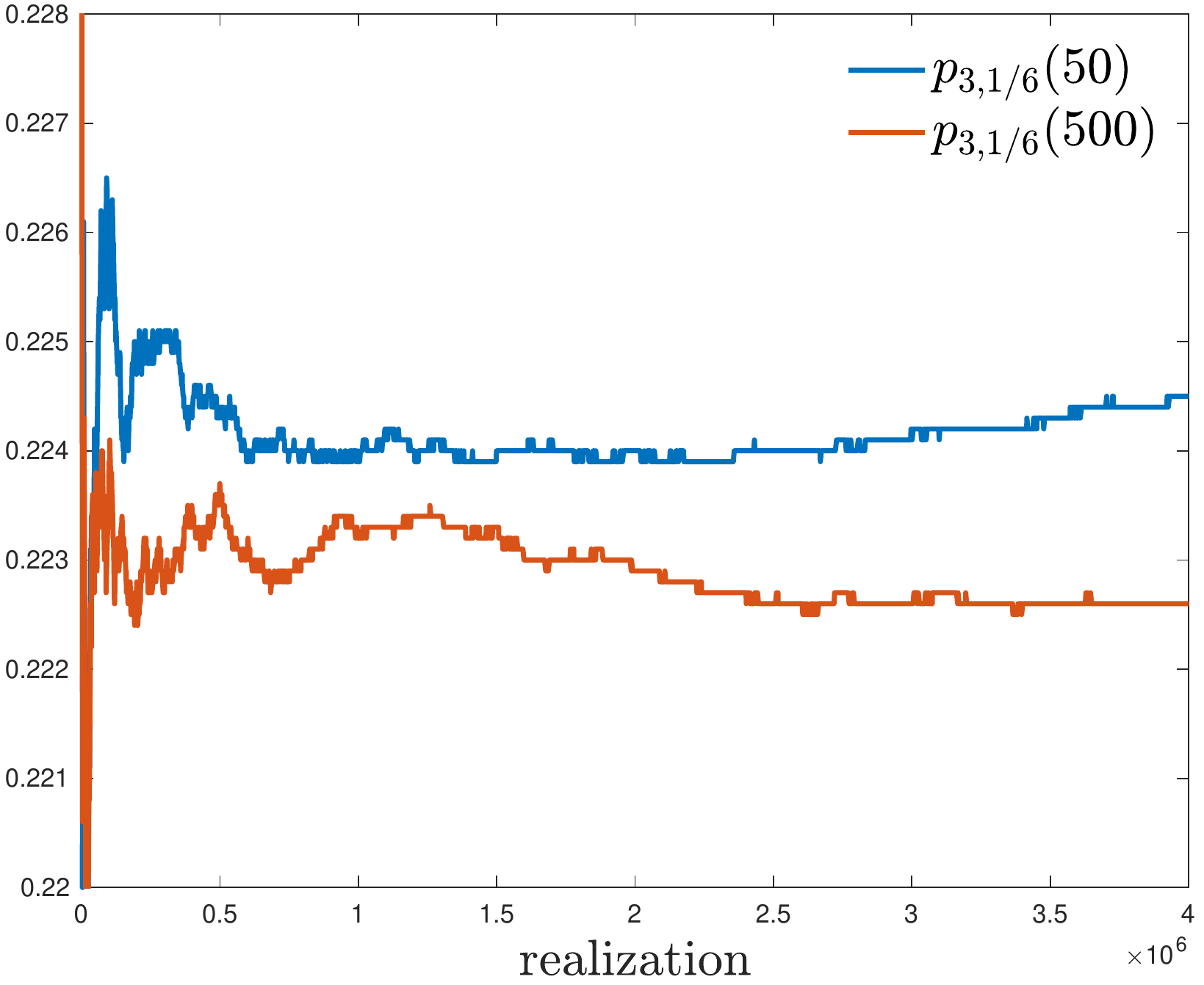}
\caption{Numerical approximations to $p_{3,1/6}(50)$ and $p_{3,1/6}(500)$ by $4\times 10^6$ i.i.d. stochastic realizations.}
\label{fig_LLN}
\end{figure}

In addition to the aforementioned theoretical proof, the low dimensional non-monotonicity can also be verified in numerical simulation. In Figure \ref{fig_LLN} we present numerical approximations of $p_{3,1/6}(50)$ and $p_{3,1/6}(500)$ achieved from $4\times 10^6$ i.i.d. stochastic realizations. 

\begin{figure}[H]
\centering
\includegraphics[width=4.8in,height=3.2in]{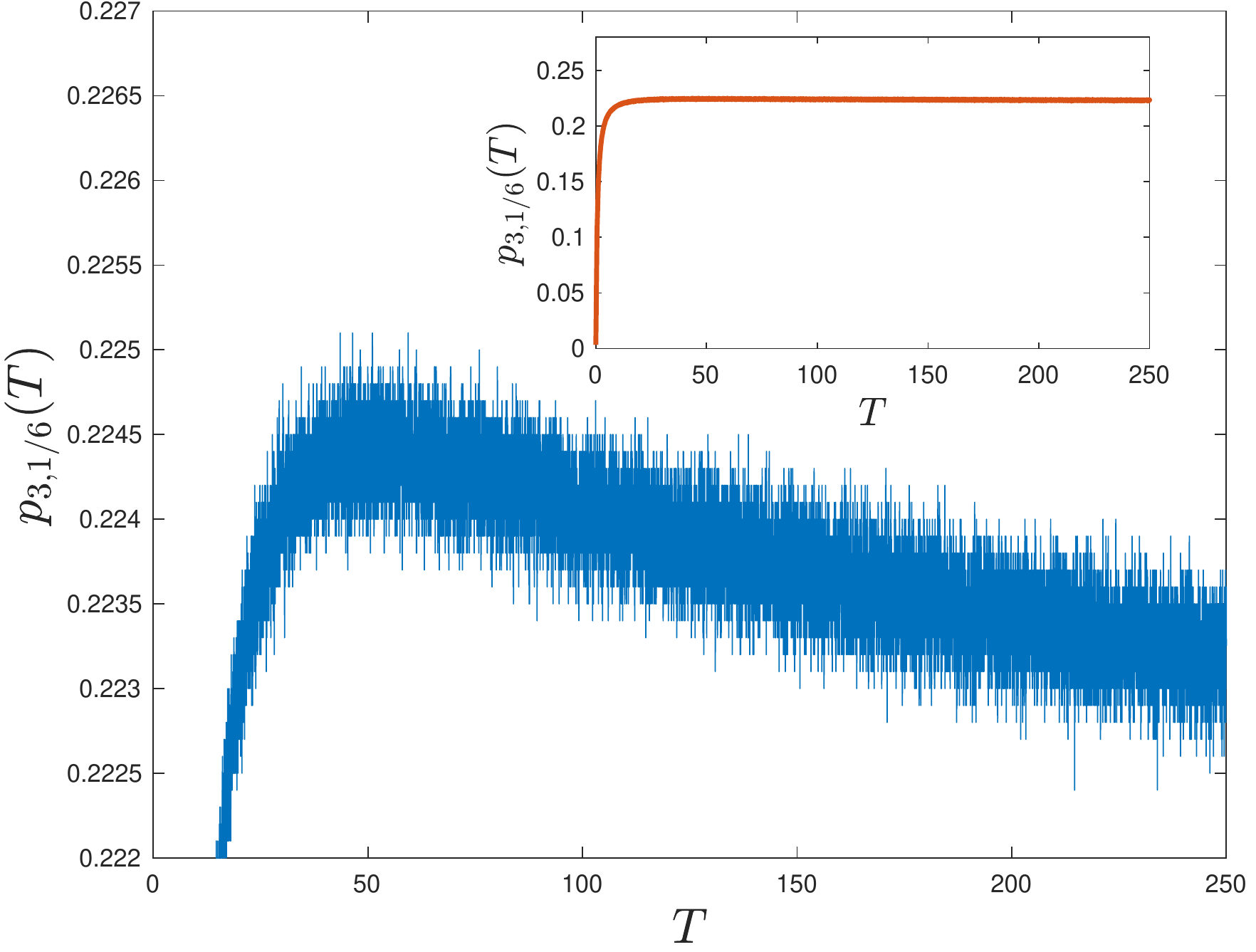}
\caption{Numerical approximations for $p_{3,1/6}(u), \ T\in [0,250]$ with spacing $\Delta T = 0.01$. Each point is evaluated by $4\times 10^6$ i.i.d. stochastic realizations. The probability seems to be non-monotonic and reaches its maximum at about $T_{\max} = 50$.  }
\label{fig_one_denisty}
\end{figure}

In Figure \ref{fig_LLN}, a significant difference between the blue and red curves is observed. Since with $4\times 10^6$ i.i.d. stochastic realizations, we have the frequencies 
$$
N_b\overset{d}{=} B\big(4\times 10^6, p_{3,1/6}(50)\big), \ N_r\overset{d}{=} B\big(4\times 10^6, p_{3,1/6}(500)\big).
$$
So their standard deviations can be bounded from above by $1/(2\times 2\times 10^3)=2.5\times 10^{-4}$. However, the difference between our approximations is about $1.7\times 10^{-3}$, which is larger than 4 times the upper bound of standard deviation. In Figure \ref{fig_one_denisty} we numerically approximate the single edge density $p_{3,1/6}(\cdot)$ for different $T$ with spacing $\Delta T = 0.01$,  and each point is evaluated by $4\times 10^6$ i.i.d. stochastic realizations. In spite of some stochastic fluctuations, the trend of non-monotonicity is clear and and the probability seems to reach maximum at $T_{max}\approx 50$.  

As dimension increases, the following theorem shows that the single edge density becomes monotonic in $T$ for all sufficiently large $d$:
\begin{Theorem}
\label{thm_d_large}
For any $u>0$, there exists $d_0=d_0(u)<\infty$ such that $p_{d,u}'(T)>0$ for all $T\in(0,\infty)$ and $d\ge d_0$. 
\end{Theorem}

The proof of Theorem \ref{thm_d_large} is similar in spirit to that of Proposition \ref{Lawler_11_1} and Theorem \ref{d34}. The main technicalities involved in the proof are asymptotic estimates for high dimensional SRW's, which are not directly related to the main scope of this paper. So we leave it in Appendix \ref{app_A}. 

Theorem \ref{thm_d_large} also motivates us to propose the following open question:

\begin{question}
For $u>0$, is there a $d_1=d_1(u)$ such that $\mathcal{FI}_d^{u,T}$ is stochastically increasing in $T$ for all $d\ge d_1$?  
\end{question}

\section{Characterization of phase diagram}
\label{sec_num}

In this section, we focus on the (potential) phase diagram of edge percolation in FRI. We start with the following result that partly characterizes the supercritical phase. 
\begin{Theorem}
\label{thm_sup_part}
For all $d\ge 3$ and FRI $\mathcal{FI}_d^{u,T}$ we have
\begin{enumerate}
\item If $\mathcal{FI}_d^{u,T}$  has an infinite cluster almost surely, then so does $\mathcal{FI}_d^{u',T}$ for all $u'>u$. 
\item (Theorem 1,  \cite{FRI_1}) For all $u >0$, there is a $0 < T_1 (u,d) < \infty$ such that for all $T > T_1$, $\mathcal{FI}_d^{u,T}$ has an infinite cluster almost surely.
\item Let $p^c_d$ be the critical edge density for $d-$dimensional Bernoulli bond percolation. For any $u>-2\log(1-p_d^c)$, there exist some $\delta=\delta(u,d)>0$ such that $\mathcal{FI}_d^{u,T}$ has an infinite cluster almost surely for $T\in [(1+\delta)^{-1}, 1+\delta]$. Moreover, for any fixed $d$, $\delta(u,d)\asymp u$ as $u\to\infty$. 

\item For any $d\ge 3$, there is $U_d<\infty$, such that for all $u\ge U_d$,  $\mathcal{FI}_d^{u,T}$ has an infinite cluster almost surely for all $T\ge (1+\delta(u,d))^{-1}$. 
\end{enumerate}
\end{Theorem}

\begin{proof}
Claim (i) is an immediately result of the monotonicity of $\mathcal{FI}_d^{u,T}$ with respect to $u$. For Claim (iii) and (iv), the key idea is to bound $\mathcal{FI}_d^{u,T}$ from below by a supercritical Bernoulli percolation. Without loss \textcolor[rgb]{1,0,0}{of} generality, one may first consider vertex $0$, edge $e_1$ and the collection of fibers with length $\ge1$ that start from $0$ and traverse $e_1$ in their first jump. We denote the number of such paths by $N_{e_1,+}$. Recalling the definition of $\mathcal{FI}_d^{u,T}$, there are $Pois(2du/(T+1))$ fibers starting from 0. While for each of them, the probability it has length at least 1 is $T/(T+1)$, and the probability it takes $e_1$ in the first step is $(2d)^{-1}$. Thus by the thinning property of Poisson distribution, we have 
$$
N_{e_1,+}\overset{d}{=}Pois(uT/(T+1)^2).
$$
Similarly, one can define $N_{e_1,-}$ to be the number of fibers that start from $x_1$ and traverse $e_1$ in their first jump. By independent increment property of PPP, $N_{e_1,-}$ is independent and identically distributed as $N_{e_1,+}$. Define event 
$$
\{\text{$e_1$ is good}\}\overset{\Delta}{=} \{N_{e_1,+}+N_{e_1,-}>0\}. 
$$
Moreover for any edge $e=\{x,y\}\in \BZ^d$, one can also define $N_{e,+}, \ N_{e,-}$ in the exact same way. And thus $\{N_{e,\pm},e\in \BZ^d\}$ form a i.i.d. sequence of $\text{Poisson}(uT/(T+1)^2)$. Again define 
$$
\{\text{$e$ is good}\}\overset{\Delta}{=} \{N_{e,+}+N_{e,-}>0\}. 
$$
Thus the collection of good edges by definition forms a Bernoulli bond percolation with single edge density
\begin{equation}
\label{eq_edge_open}
p=1-P(N_{e_1,+}=0)\cdot P(N_{e_1,-}=0)=1-\exp(-2uT/(T+1)^2),
\end{equation}
which percolates when
$$
\frac{uT}{(T+1)^2}\ge -\log(1-p_d^c)/2. 
$$
Note that a good edge is by definition always traversed by the FRI. Claim (iii) is now a direct result of \eqref{eq_edge_open}, the fact that $T/(T+1)^2$ reaches its maximum of $1/4$ at $T=1$, and that $T/(T+1)^2\asymp T^{-1}$ as $T\to\infty$. 

Now for (iv), note that for a fixed $u$, say $u=1$, by Theorem 1,  \cite{FRI_1}, there is a $T_1$ such that for all $T > T_1$, $\mathcal{FI}_d^{1,T}$ has an infinite cluster almost surely. With Claim (i) we now know this also holds for all $u\ge1$. In (iii), we have shown that $\delta(u,d)\asymp u$. Thus there is always a $U_d$ such that $\delta(u,d)\ge T_1$ for all $u\ge U_d$. And thus, we have an  infinite cluster almost surely for all $T$ from $(1+\delta(u,d))^{-1}$ all the way to infinity. 
\end{proof}

In particular, when $d=3$, $u=1$, the following result is a direct corollary of Claim (iii), Theorem \ref{thm_sup_part} and the fact that $p_3^c\approx 0.25$:

\begin{Corollary}
	$\mathcal{FI}_3^{1,T}$ has a infinite cluster almost surely for all $T\in [1/4, 4]$. 
\end{Corollary}

We can replace the critical value $p^c_d$ with $p_d^{fin}$, the critical value for non-existence of infinite connected set in the complement of the percolation cluster of a Bernoulli bond percolation (see \cite{grimmettpercolation} for details), and the infinite cluster in Claim (iii) and (iv) becomes unique:

\begin{Corollary}
	For all $d\ge 3$ and FRI $\mathcal{FI}_d^{u,T}$ we have
	\begin{enumerate}
		\item For any $u>-2\log(1-p_d^{fin})$, there exist some $\hat\delta=\hat\delta(u,d)>0$ such that $\mathcal{FI}_d^{u,T}$ has a unique infinite cluster almost surely for $T\in [(1+\hat\delta)^{-1}, 1+\hat\delta]$. Moreover, for any fixed $d$, $\hat\delta(u,d)\asymp u$ as $u\to\infty$. 
		
		\item For any $d\ge 3$, there is $\hat U_d<\infty$, such that for all $u\ge \hat U_d$,  $\mathcal{FI}_d^{u,T}$ has a unique infinite cluster almost surely for all $T\ge (1+\hat\delta(u,d))^{-1}$. 
	\end{enumerate}
\end{Corollary}

For $u\in (-2\log(1-p_d^c), U_d)$, there appears to be a gap between the known supercritical regime in Claim (iii) and Theorem 1,  \cite{FRI_1}. And in the moderate $T$ regime, the percolation is more "Bernoulli like", while in the large $T$ regime, the percolation is more "interlacement like". In Theorem \ref{d34}, it has been proved that low dimensional FRI does not enjoy stochastic monotonicity. This casts more doubt on whether the transition between the subcritical phase for very small $T$ and the supercritical phase for sufficiently large $T$ is sharp or unique. In the rest of this section, we make numerical explorations towards this direction. The general guidelines behind criteria of numerical tests in this section are mostly inspired by \cite{MR1707339,MR2838338}.

In order to investigate the phase transition property of FRI, we develop the following parallel computing algorithm in order to efficiently sample the configuration within a large box in $\BZ^d$ of size $N$, with data transferred within up to $80$ cores via the Message Passing Interface (MPI). In Section 4.1 of \cite{FRI_1}, it has been shown that we can sample the configuration of FRI restricted within an infinite set $K$ in the following steps:

\begin{algorithm}
\caption{Finitary Random Interlacement}

\begin{enumerate}
	\item Divide the vertices $x \in V$ into $N_p$ mutually independent batches $(B_1(N), \dots, B_{N_p}(N))$ and distribute one batch to one processor. 

	\vspace{2mm}

	\item For $s$-th batch, for any $x\in B_s(N)$, sample an independent random variable $N_x\sim Pois(\frac{2du}{T+1})$. Then sample a sequence of i.i.d. random walks $\{\eta_i\}_{i=1}^{N_x}$ independent to $N_x$, with distribution $P_x^{(T)}$.

	\vspace{2mm}

	\item For each trajectory $\eta_i, \ i\le N_x$ of the random walk mentioned above, if $\eta$ escapes from $K$ (i.e. for any $n\ge 1$, $\eta(n)\notin K$), then start a new independent random walk with distribution $P_x^{(T)}$ and collect its trajectory $\hat{\eta}_i$; if not, jump to the next trajectory $\eta_{i+1}$. 

	\vspace{2mm}

	\item Collect all the trajectories $\bigcup\limits_{s, i}\hat{\eta}_i\cap K$ from all processors. \end{enumerate} 
\end{algorithm}

%\note{Yunfeng: would you want to add some description of the code here?}

%\note{Note that the $u$ in this paper should be the $u$ in our code divided by $2d$.}

\begin{figure}[!h]
\centering
\subfigure[u = 0.1]{{\includegraphics[width=0.49\textwidth,height=0.30\textwidth]{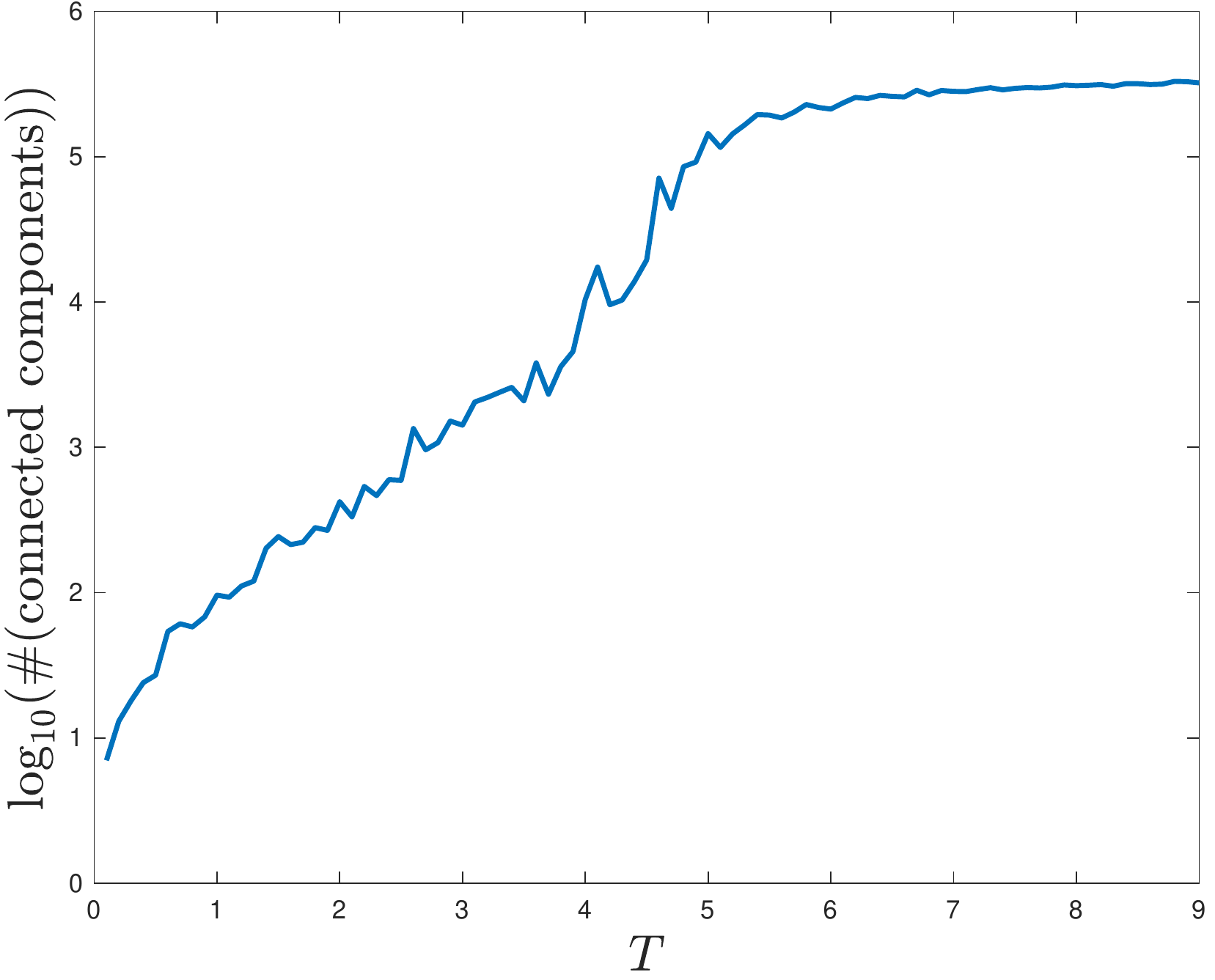}}
\includegraphics[width=0.49\textwidth,height=0.30\textwidth]{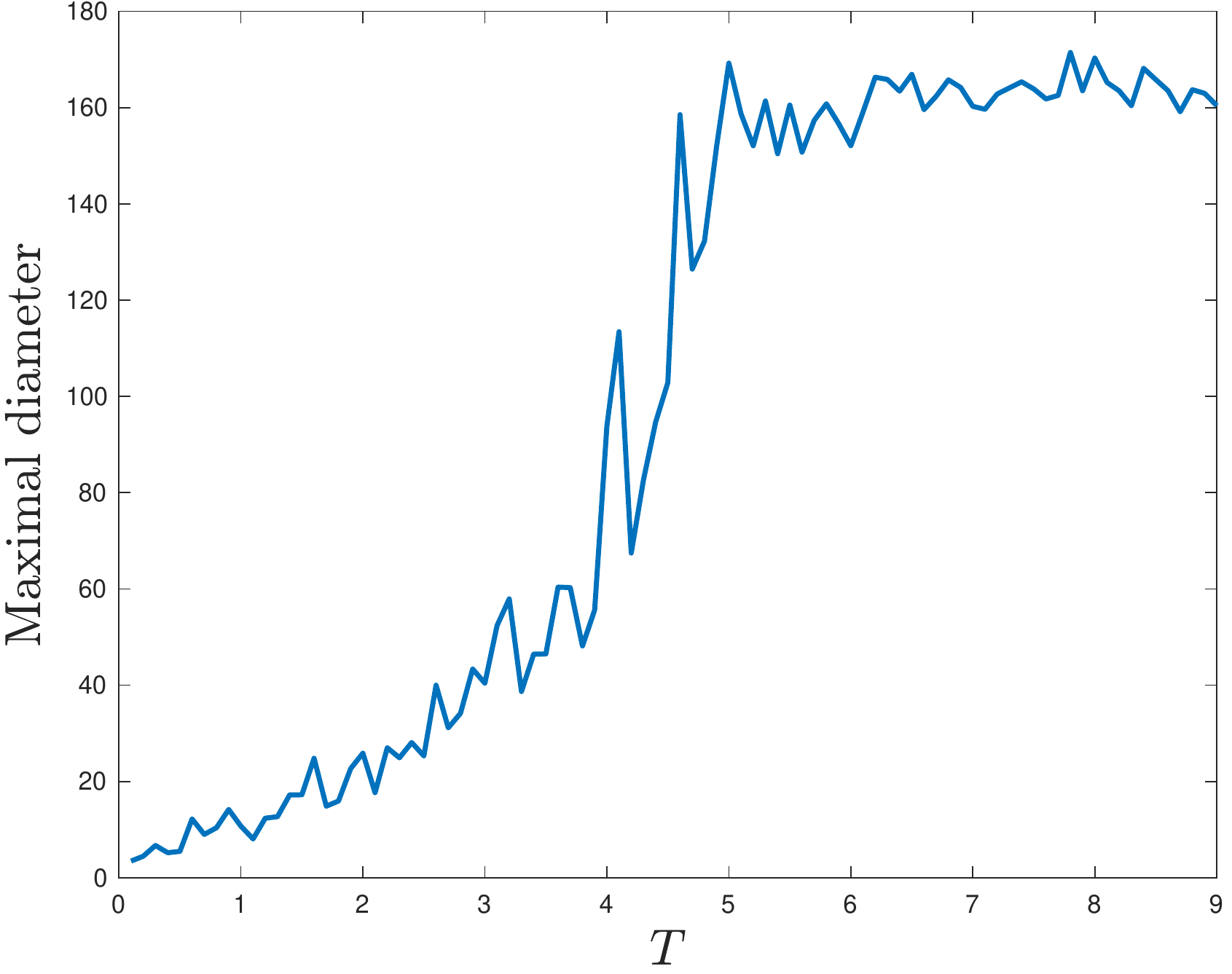}}\\
\subfigure[u = 0.2]{{\includegraphics[width=0.49\textwidth,height=0.30\textwidth]{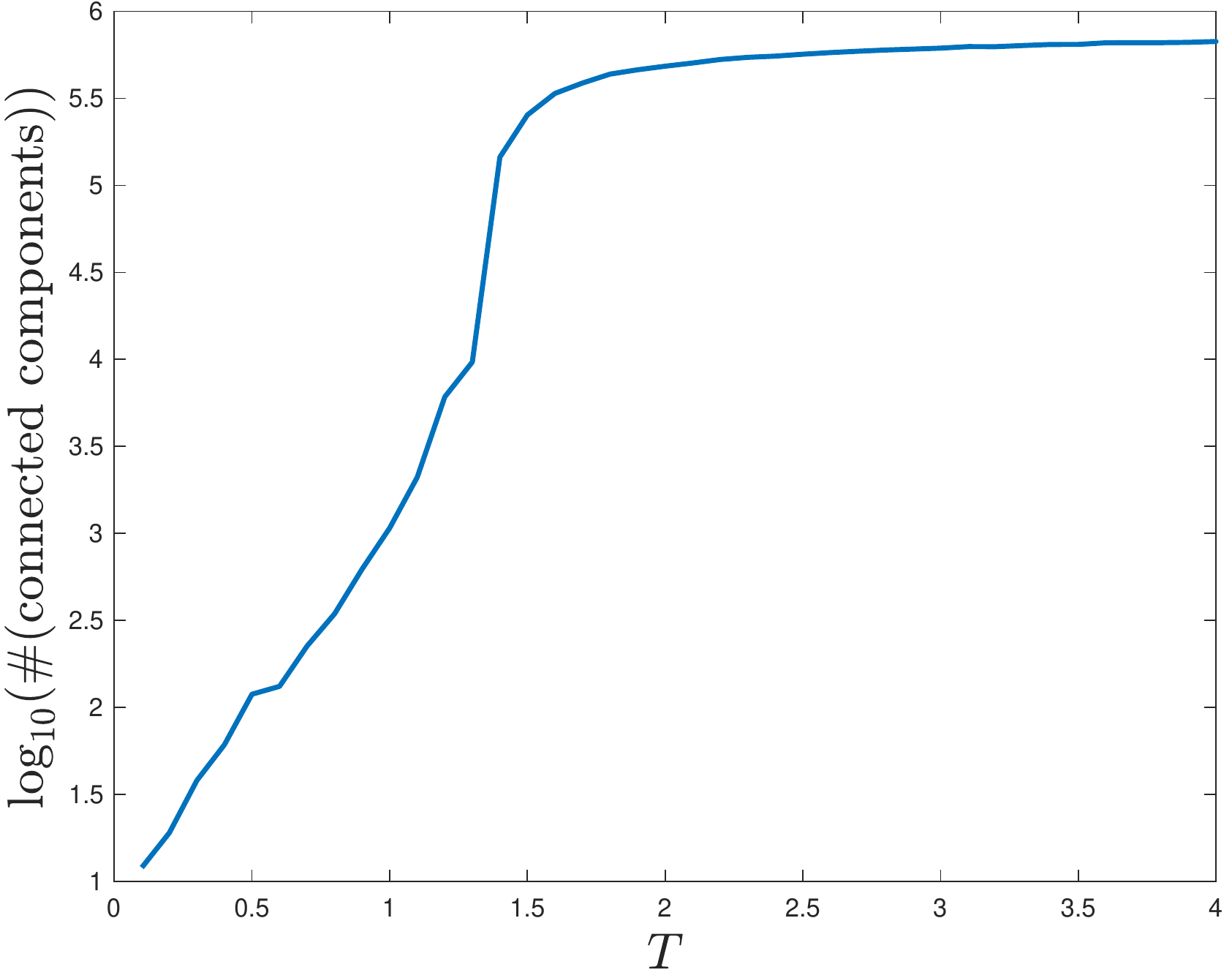}}
\includegraphics[width=0.49\textwidth,height=0.30\textwidth]{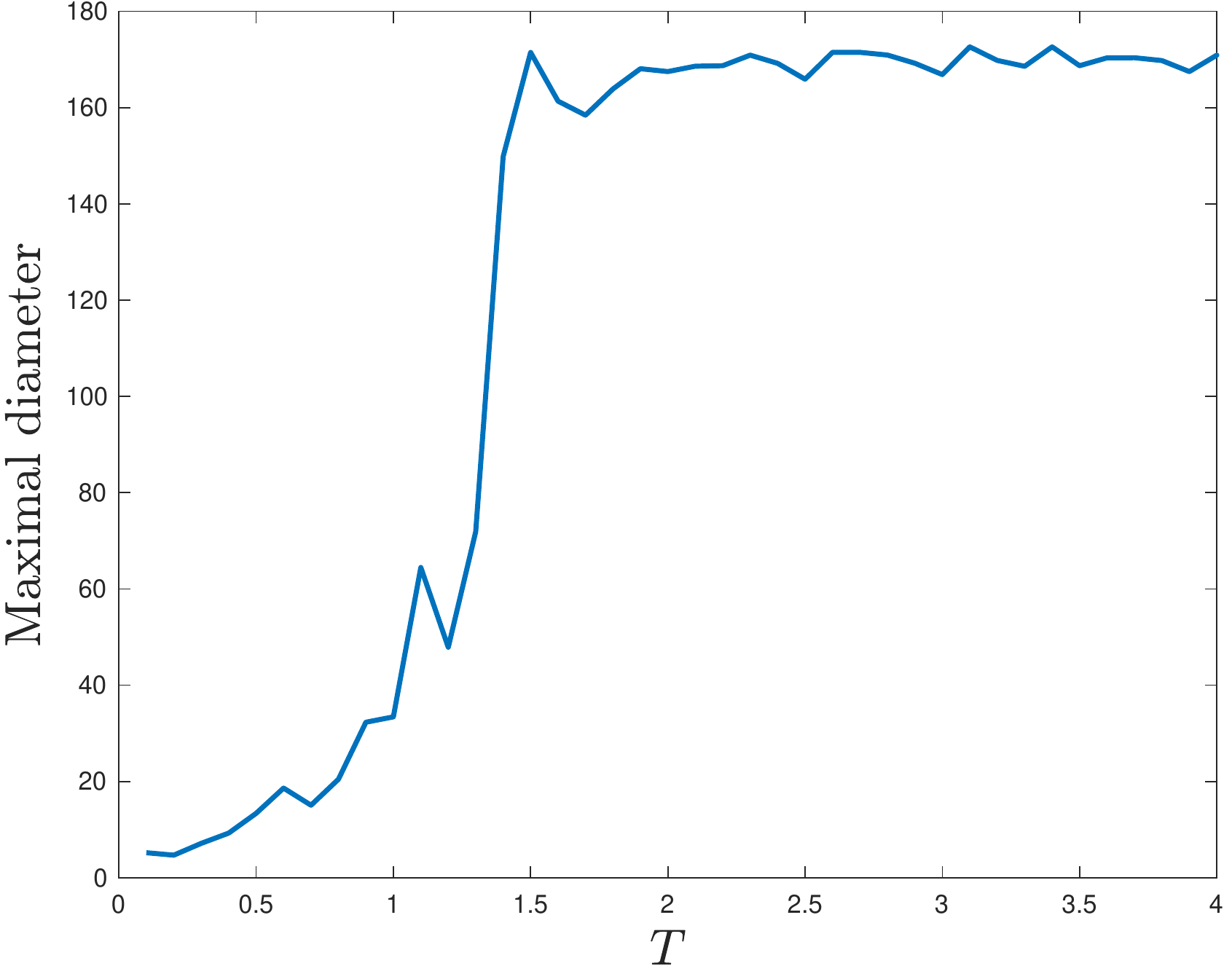}}\\
\subfigure[u = 0.5]{{\includegraphics[width=0.49\textwidth,height=0.30\textwidth]{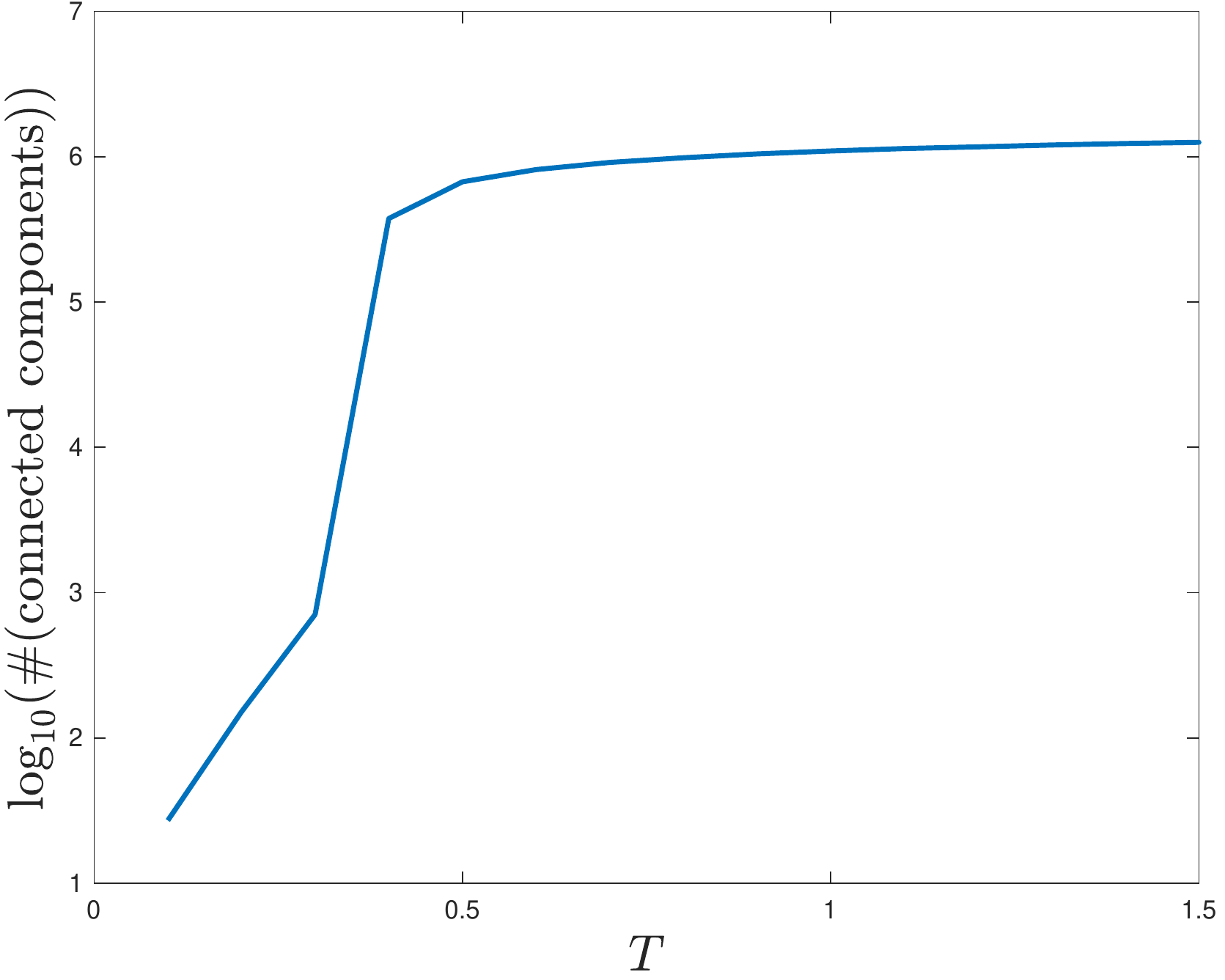}}
\includegraphics[width=0.49\textwidth,height=0.30\textwidth]{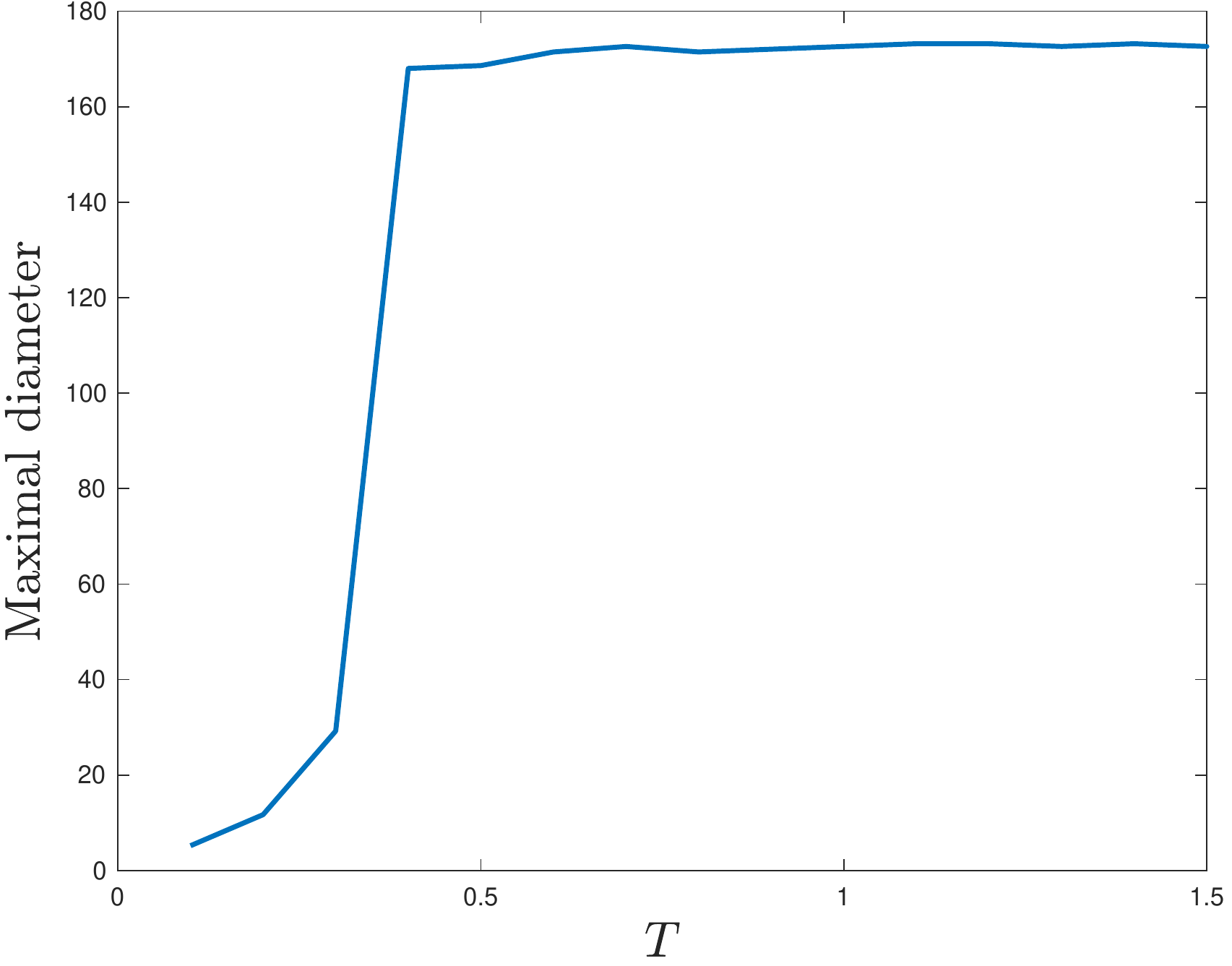}}\\
\caption{Stochastic simulations on the sizes (left, in logarithm scale) and maximal diameters (right) of FRI largest connected components. Phase transitions are clearly observed for different $u$.}
\label{fig_sharp}
\end{figure}

Since $\bigcup\limits_{s, i}\hat{\eta}_i\cap K$ is identically distributed as the collections of all fibers in $\mathcal{FI}_d^{u,T}$ which traverse $K$, we have that $\bigcup\limits_{s, i}\hat{\eta}_i\cap K\overset{d}{=} \mathcal{FI}_d^{u,T}\cap K$. Using the aforementioned algorithm, we first explore the sharpness and uniqueness of the phase transition. For this propose, one can naturally look at the size (in either cardinality or diameter) of the largest connected component within a large box in $\BZ^d$ of size $N$, say $[0,N]^d\cap \BZ^d$. In the supercritical phase, there should be a macroscopic largest connected component within $[0,N]^d\cap \BZ^d$, since it should, with high probability, be the largest cluster in the intersection(s) between an infinite cluster and the box. While in the subcritical phase, the largest connected component should be microscopic with respect to $N$. See Figure \ref{fig_1a}  for illustration. In Figure \ref{fig_sharp}, we present stochastic simulation results on the cardinality and maximal diameter of the FRI's largest connected components within  $[0,N]^3\cap \BZ^3$, for $N=100$, $u=0.1,0.2,0.5$ and various $T$'s, under only one realization. In order to manifest the phase transition more clearly, we choose different ranges of $T$ under different $u$. Although the curves are not smooth due to some random fluctuations and size effects, the phase transitions are clearly observed, and the critical value $T_c$ seems to be smaller as $u$ goes larger.

%\note{6 subfigures, 2-by-3, Subfig. 1.(a) for $u=0.1$ cardinality (maybe log), Subfig. 1.(b) for $u=0.1$ diameter, etc. }

\begin{figure}[!h]
\centering
\subfigure[Cardinality of connected components (in log10 sacle)]{\includegraphics[width=0.49\textwidth,height=0.30\textwidth]{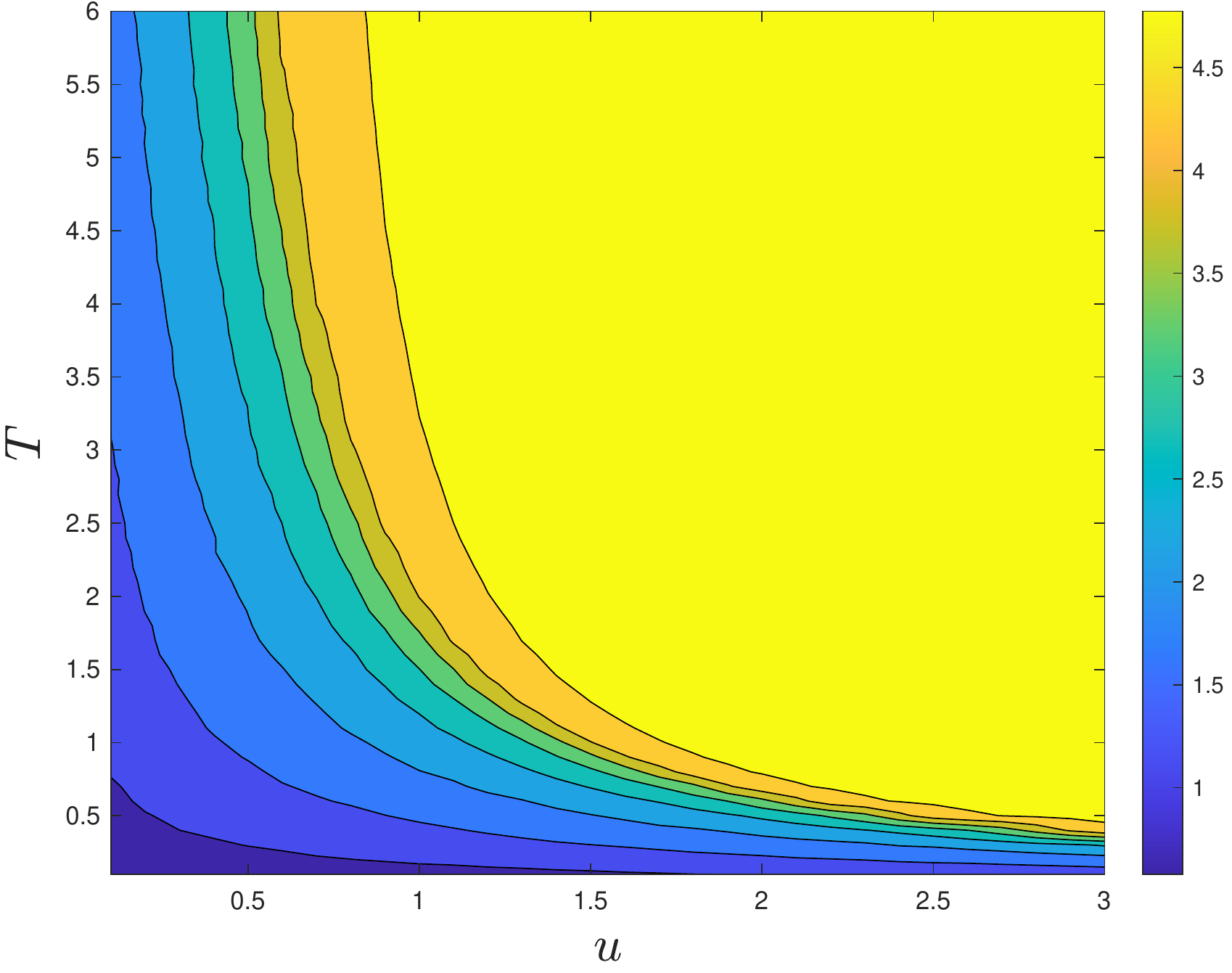}}
\subfigure[Maximal diameter]{\includegraphics[width=0.49\textwidth,height=0.30\textwidth]{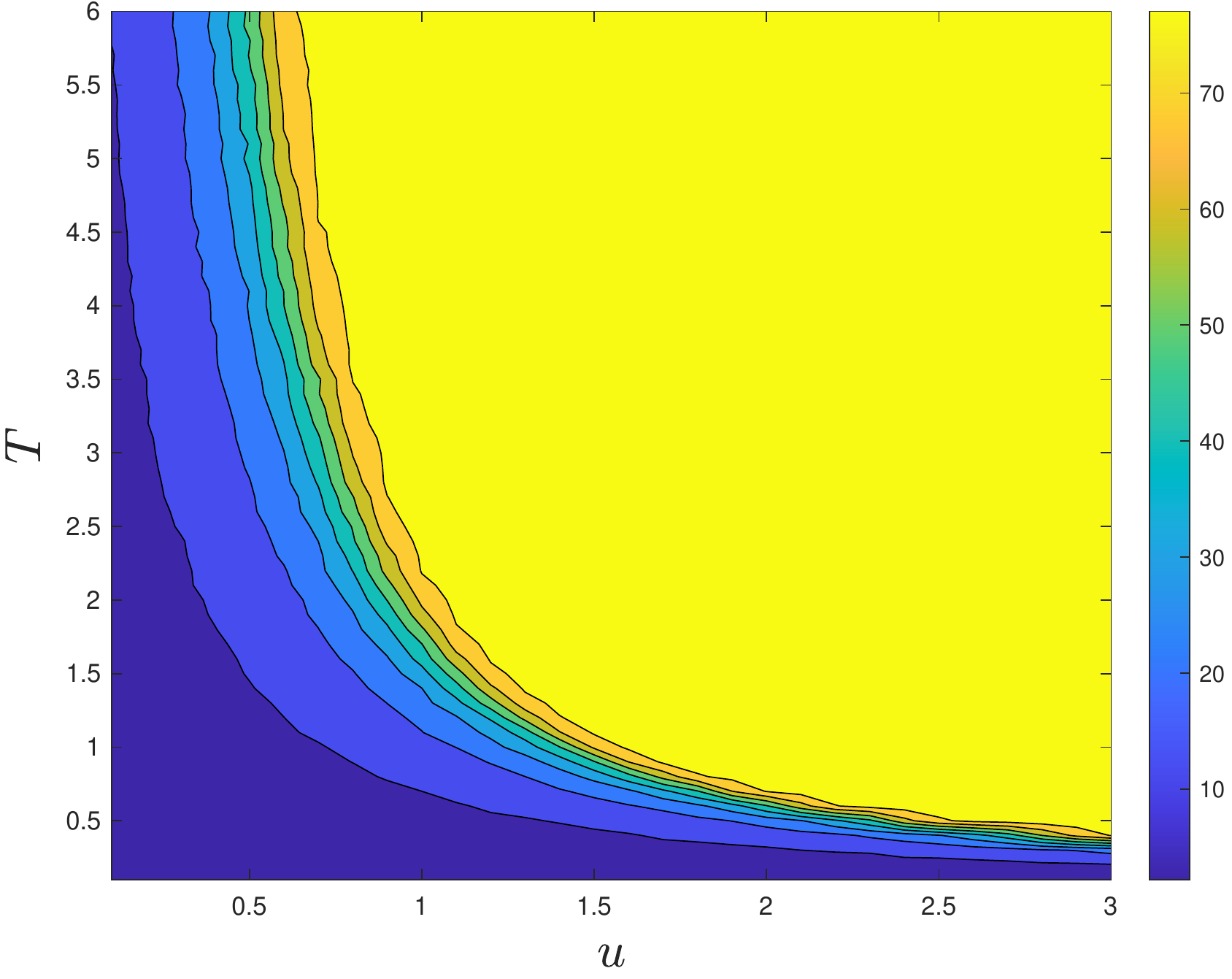}}
\caption{Illustration for the existence of phase diagram of FRI largest connected components. Each point is obtained by the averaging of 100 i.i.d FRIs. The sharpness and uniqueness of phase transition seem to hold for all $u$'s.}
\label{fig_diag}
\end{figure}

Figure  \ref{fig_sharp} strongly supports that there exists a unique and sharp phase transition in $T$ for the existence of infinite connected cluster(s), at least for the $u$'s we chose. This encourages us to extend the test for all combinations of $(u,T)$'s within an appropriate grid $[0, 3] \times [0, 6]$ with $\Delta u=0.1$, $\Delta T = 0.1$. Noting that the computational cost grows as $\mathcal{O}(N^6)$, we need to work on a smaller $N=50$. And in order to avoid the extra randomness due to the smaller box size, we run 100 i.i.d. FRI copies for each combination of $(u,T)$ and approximate the expected size of the largest cluster. Results shown in Figure \ref{fig_diag} indicate that the sharpness and uniqueness of phase transition seem to hold for all $u$'s.

%\begin{figure}[!h]
%\centering
%\subfigure[Number of connected components (in log10 sacle)]{\includegraphics[width=0.49\textwidth,height=0.30\textwidth]{1st_ncc_phase_2d_old.pdf}}
%\subfigure[Maximal diameter]{\includegraphics[width=0.49\textwidth,height=0.30\textwidth]{1st_max_diameter_phase_2d_old.pdf}}
%\caption{Illustration for the existence of phase diagram}
%\label{fig_diag}
%\end{figure}
%\note{2 subfigures one for log-cardinality one for diameter}

\begin{figure}[!h]
\centering
\subfigure[u = 0.1]{{\includegraphics[width=0.49\textwidth,height=0.30\textwidth]{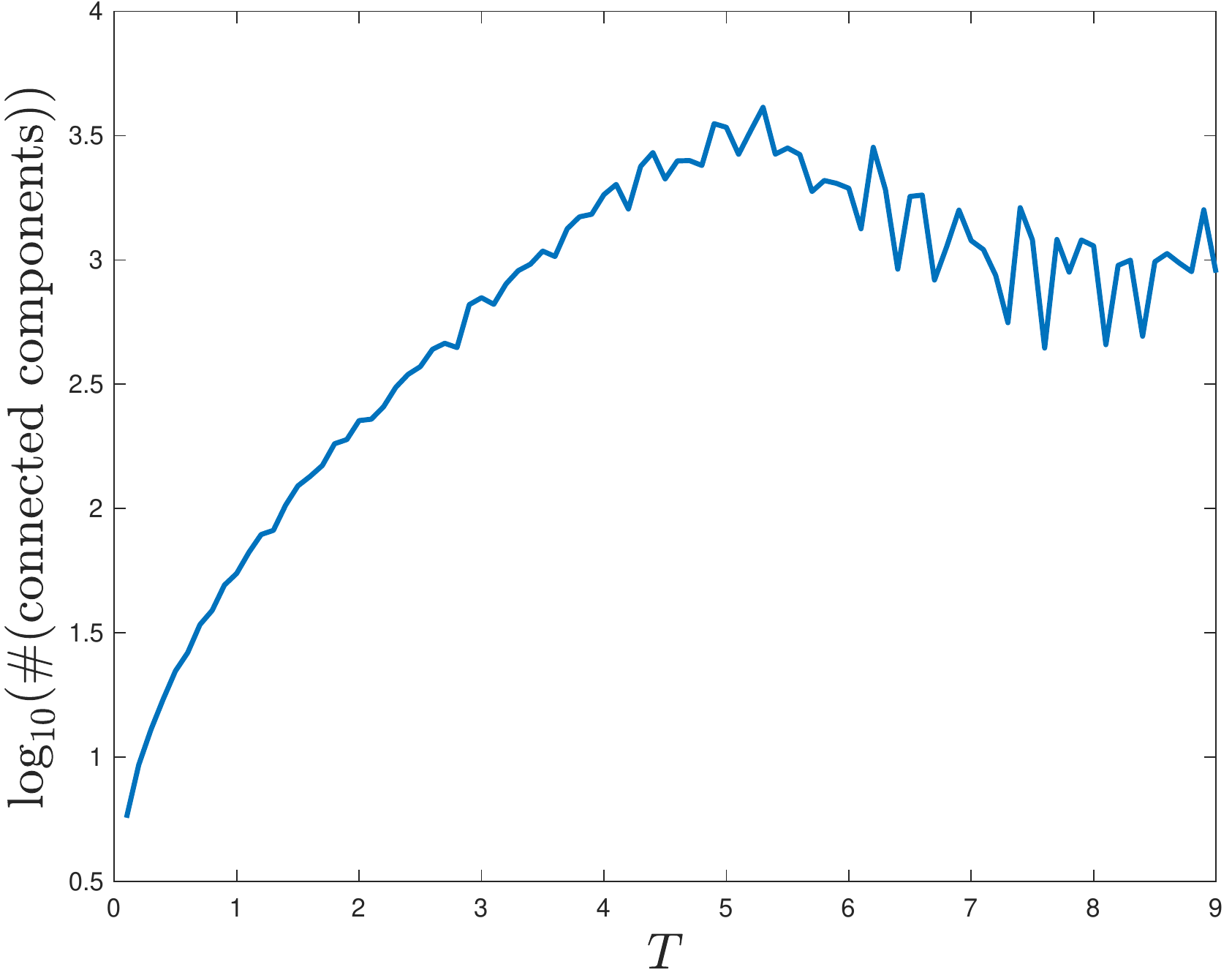}}
\includegraphics[width=0.49\textwidth,height=0.30\textwidth]{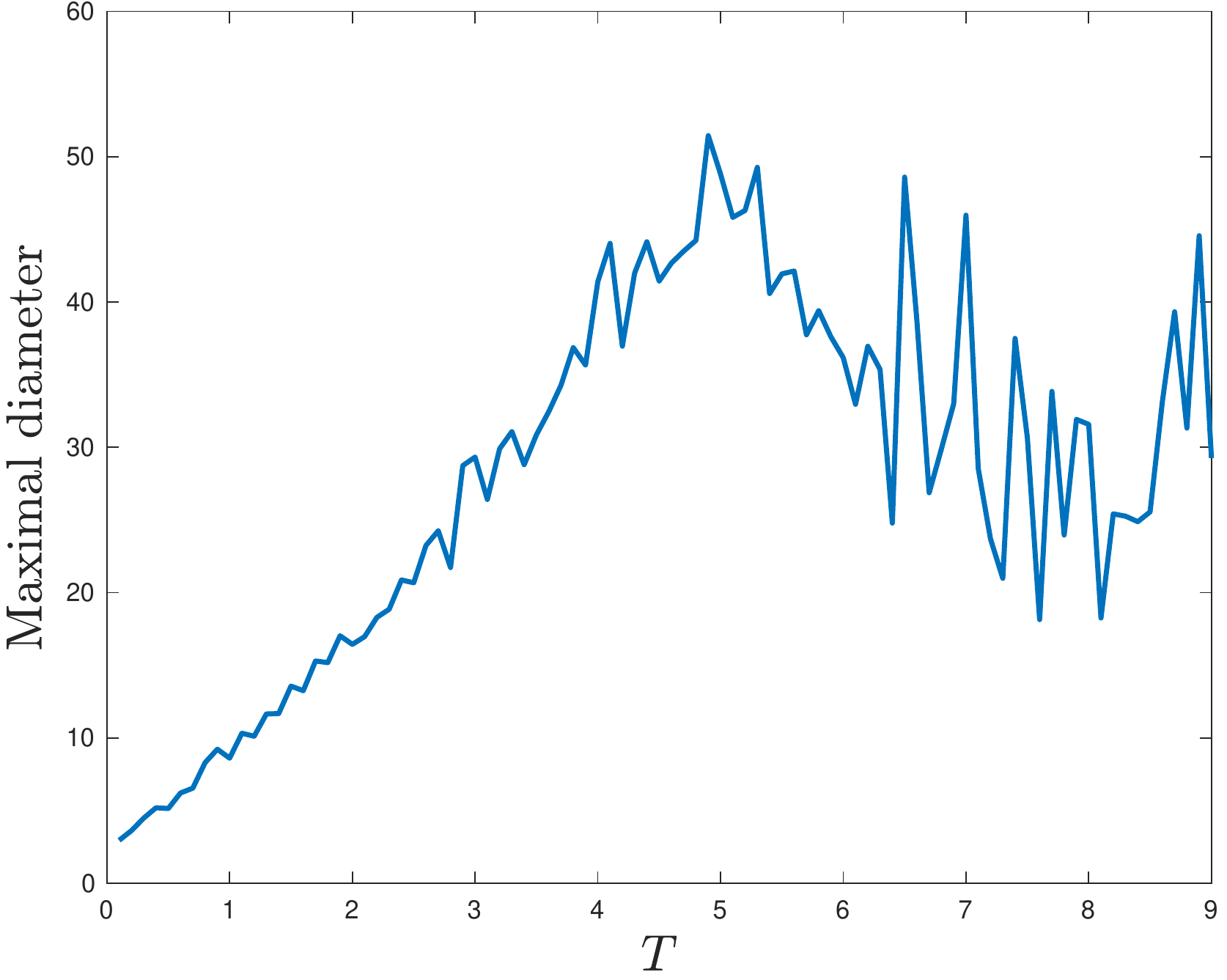}}\\
\subfigure[u = 0.2]{{\includegraphics[width=0.49\textwidth,height=0.30\textwidth]{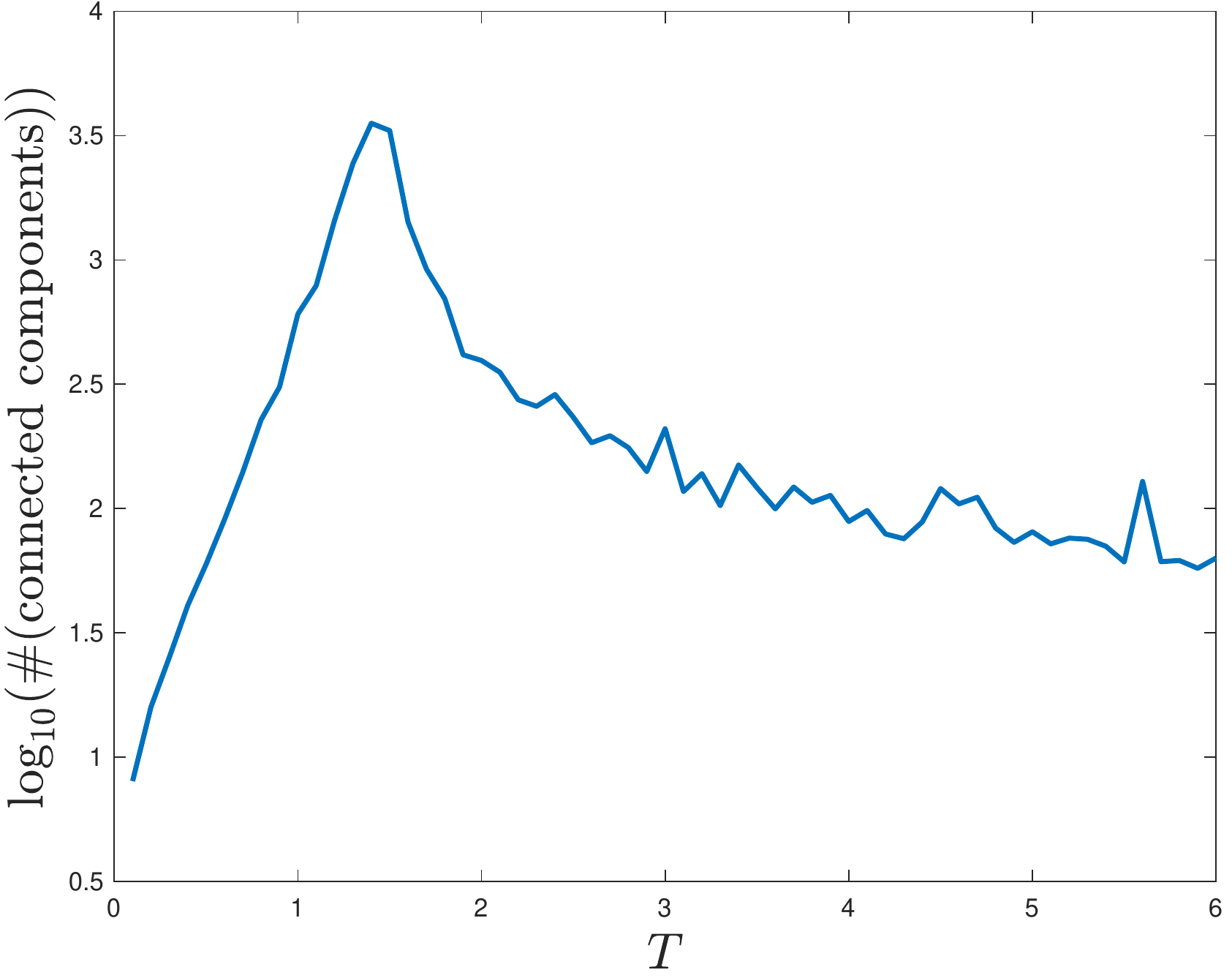}}
\includegraphics[width=0.49\textwidth,height=0.30\textwidth]{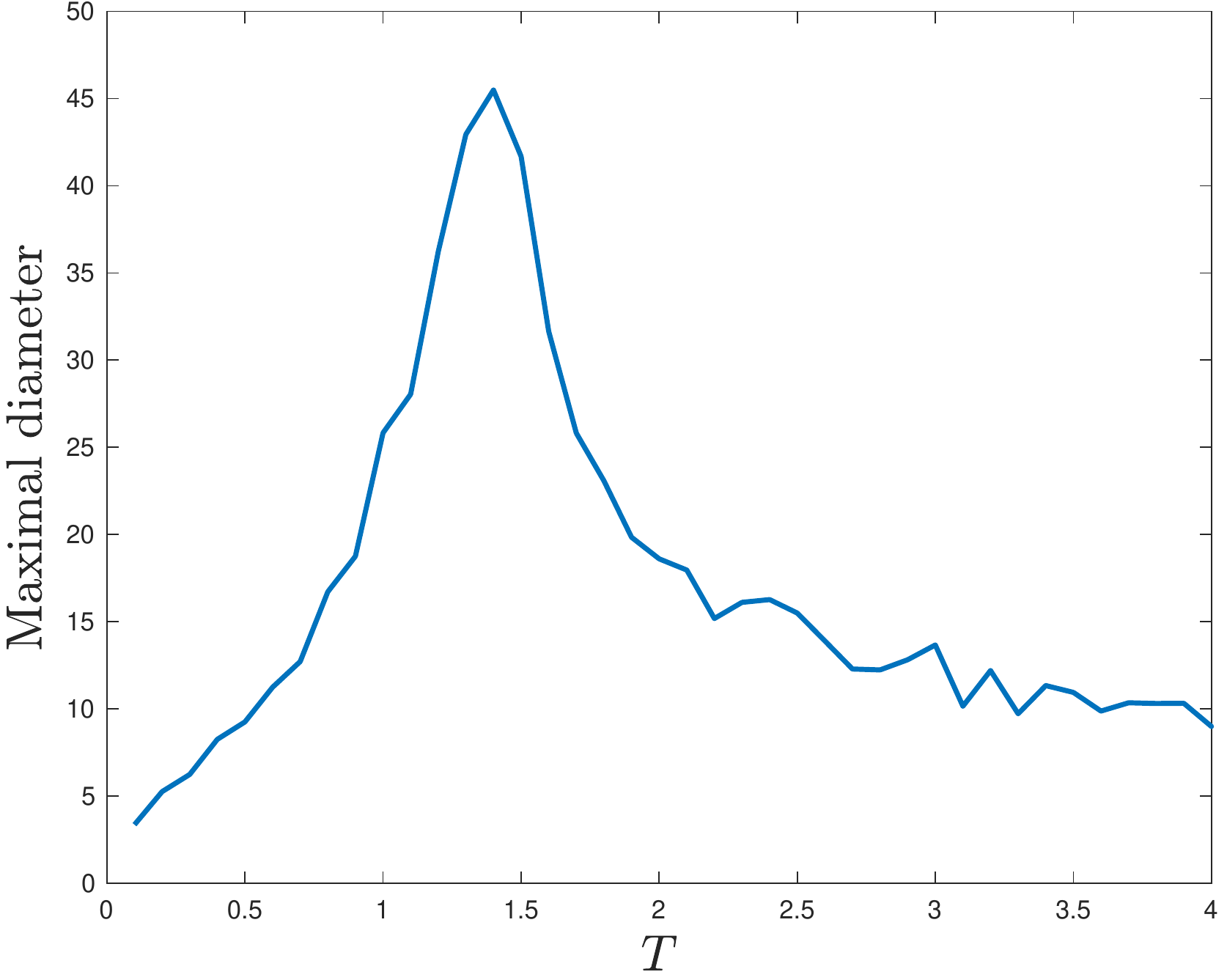}}\\
\subfigure[u = 0.5]{{\includegraphics[width=0.49\textwidth,height=0.30\textwidth]{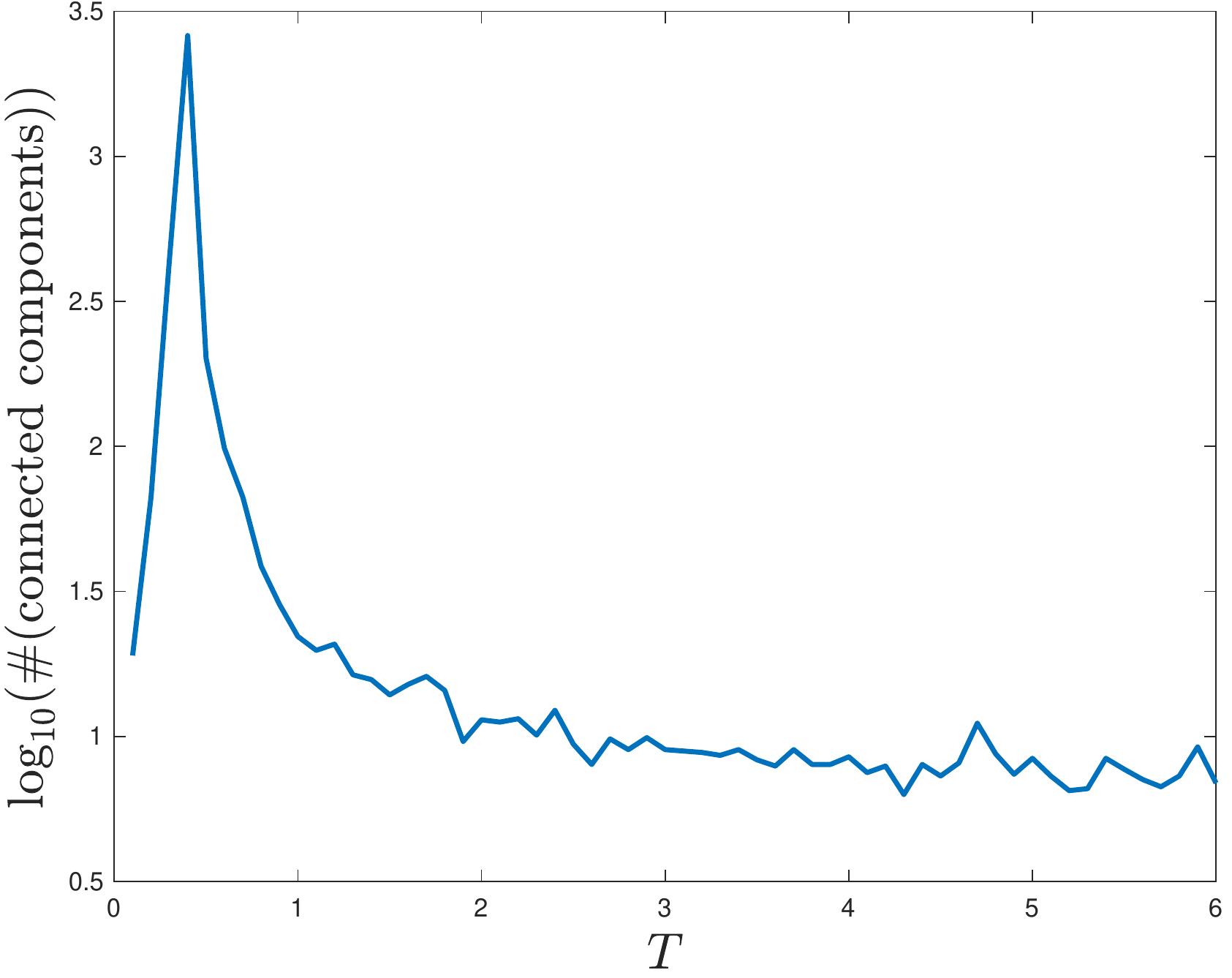}}
\includegraphics[width=0.49\textwidth,height=0.30\textwidth]{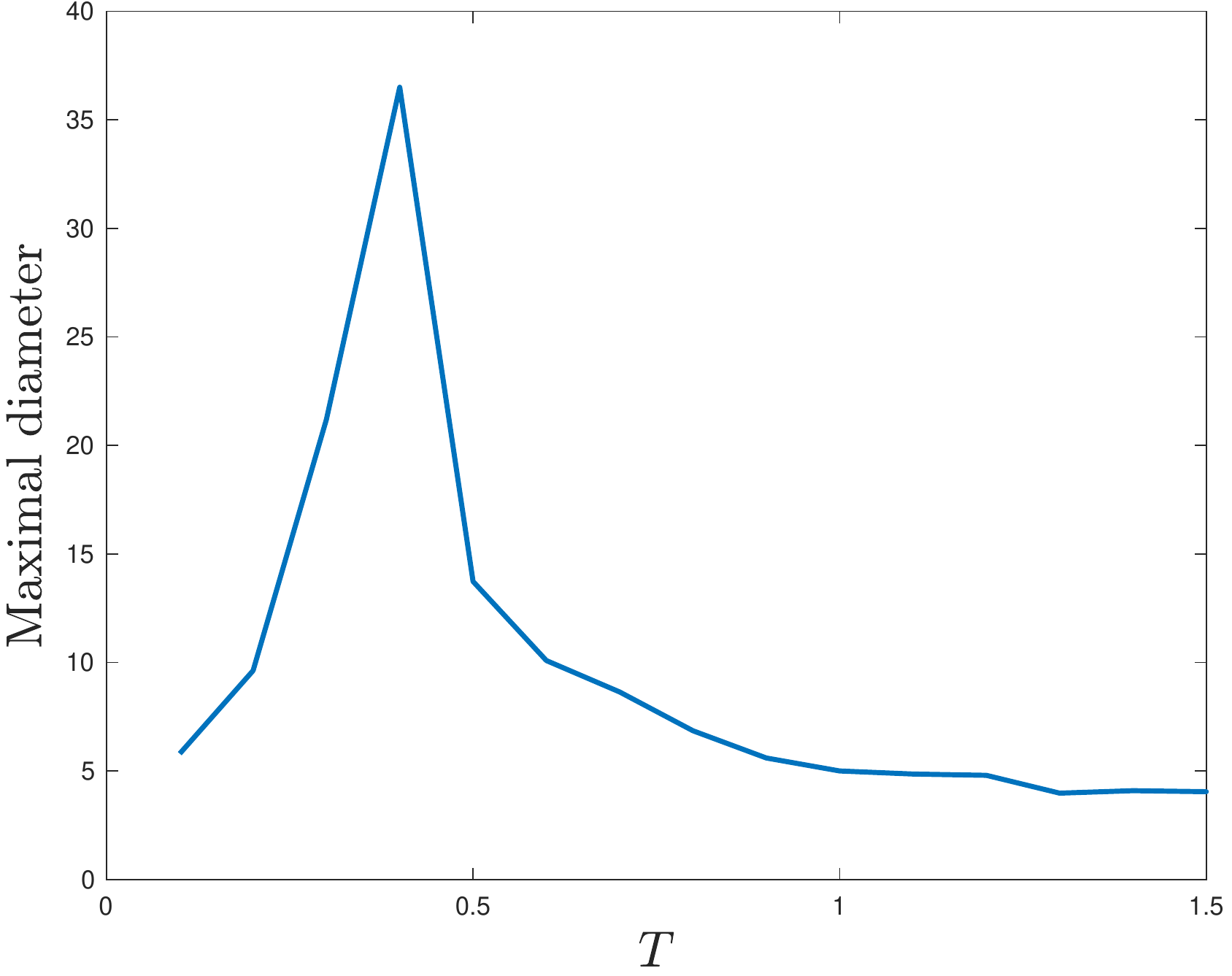}}\\
\caption{Stochastic simulations on the sizes (left, in logarithm scale) and maximal diameters (right) of FRI second largest connected components. Only a sharp peak is observed at criticality, indicating that the infinite cluster in the supercritical phase is very like to be unique. }
\label{fig_second}
\end{figure}

\begin{figure}[!h]
\centering
\subfigure[Number of connected components (in log10 sacle)]{\includegraphics[width=0.49\textwidth,height=0.30\textwidth]{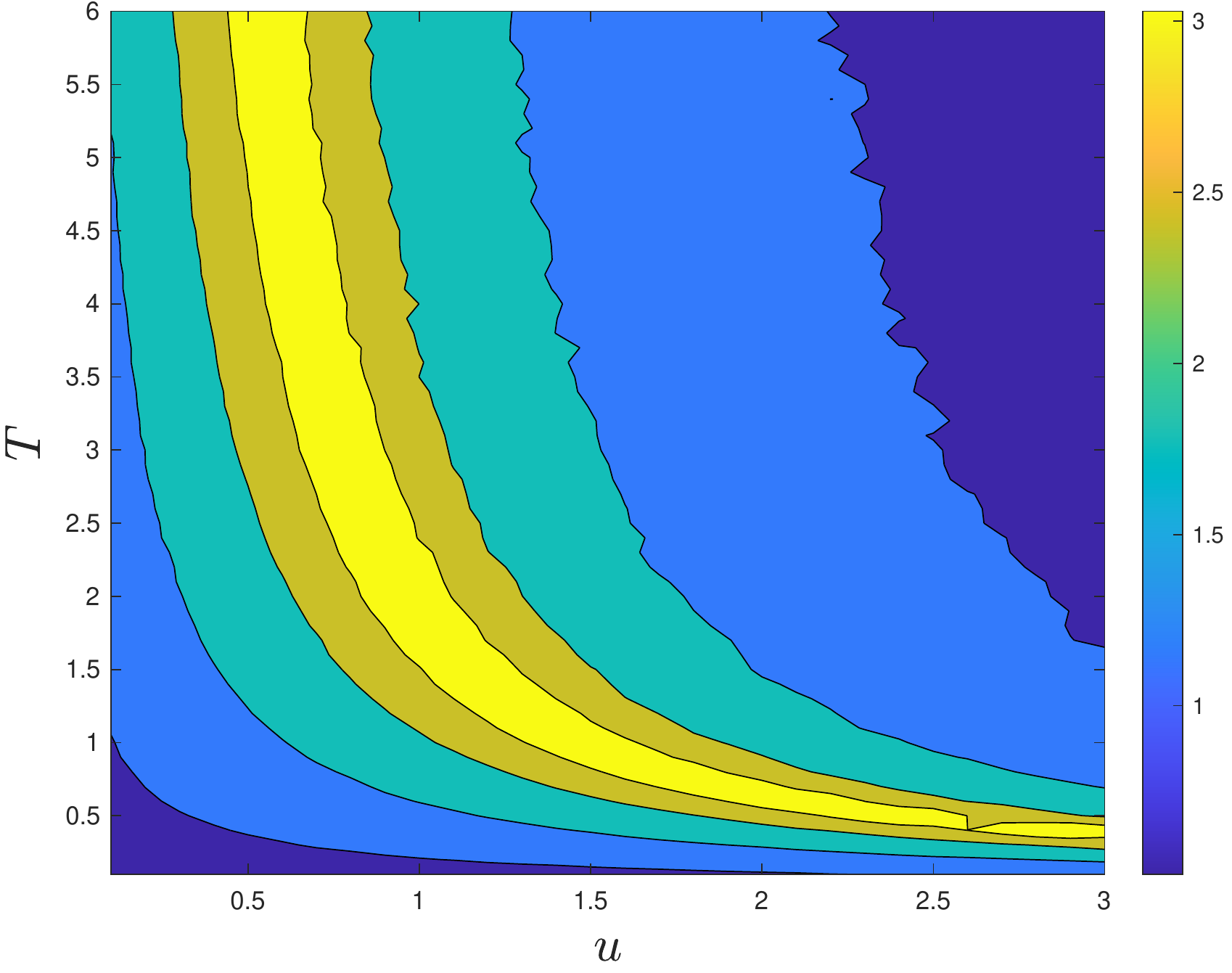}}
\subfigure[Maximal diameter]{\includegraphics[width=0.49\textwidth,height=0.30\textwidth]{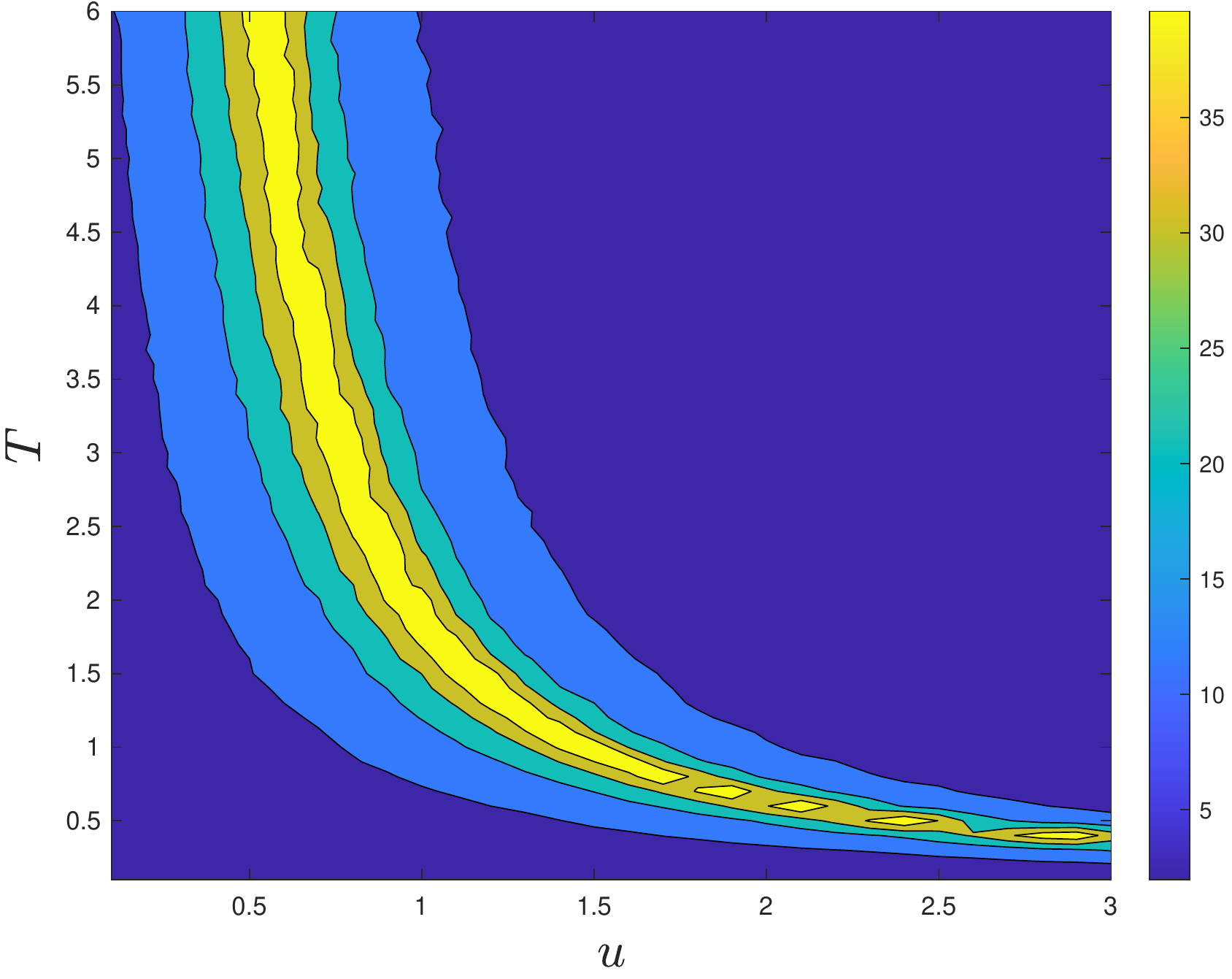}}
\caption{Illustration for the existence of phase diagram of FRI second largest connected components. Each point is obtained by the averaging of 100 i.i.d FRIs.  Only a single band is observed at criticality, indicating that the infinite cluster in the supercritical phase is very like to be unique.}
\label{fig_diag_2nd}
\end{figure}

Now we test the uniqueness of the infinite cluster(s) in the supercritical phase by investigating the size of second largest cluster within our box. Suppose there exists a region in the phase space, where there are more than one infinite clusters. Then both clusters should intersect a sufficiently large box. Thus we should expect a macroscopic second largest cluster. However, as shown in Figure \ref{fig_second}, the second largest clusters remain universally microscopic in the phase space, especially in the supercritical phase, except for a sharp peak at criticality.  Besides, as presented in Figure \ref{fig_diag_2nd}, a single band structure is observed in the phase diagram for the second largest cluster. These strongly suggests that the infinite cluster in the supercritical phase is very likely to be unique.

%\note{2 subfigures one for log-cardinality one for diameter}

\begin{Remark}
For the peak near $T_c$, one may think that when $T$ is slightly over criticality, FRI does have an infinite cluster but a very "slim" one. And since our box is finite, it may have two "branches" intersecting our box but are only connected outsides it. 
\end{Remark}

Based on the aforementioned numerical evidences and the fact that the case $d=3$ shows the strongest non-monotonicity in the single edge density estimation, we propose the following conjecture:
\begin{conjecture}
\label{conjecture_cr}
For all $d\ge 3$ and $u>0$, there is a $T_c=T_c(u)\in (0,\infty)$ such that for $\mathcal{FI}^{u,T}_d$
\begin{itemize}
\item There is a.s. no infinite cluster for all $T<T_c$. 
\item There is a.s. a unique infinite cluster for all $T>T_c$. 
\end{itemize}
\end{conjecture}

If we for now accept the existence of critical value in the conjecture above, we conclude this section by exploring the shape and asymptotic of the curve of $T_c$. With the help of monotonicity over $u$ (Claim (i) Theorem \ref{thm_sup_part}), we use the following hill climbing algorithm and record the ascending path $(u_0, T_0), \dots, (u_n, T_n), \dots $, with small spacing $\Delta T = 0.01$, $\Delta u = 0.01$, $u_0 = 3$, $T_0 = 0.01$. This algorithm significantly reduces the numerical costs in finding the boundary of phase transitions.

\begin{algorithm}[H]
\caption{Hill climbing algorithm\label{hill_climbing}}

\

\textbf{Input:} The box size $N$, A sufficiently small initial $T_0$, a sufficiently large initial $u_0$, spacings $\Delta T$ and $\Delta u$ and a threshold $\varepsilon$.

\begin{enumerate}

	\item Start from $n= 0$ and $(u_0, T_0)$.

	\vspace{2mm}
	
	\item Run one realization of FRI under the parameters $(u_n, T_n)$ and calculate the maximal diameter $d_n$ of the largest cluster of such FRI.

	\vspace{2mm}
	
	\item If $d_n < \sqrt{3}\varepsilon N$, $(u_{n+1}, T_{n+1}) = (u_{n}, T_{n} + \Delta T)$, go to Step 5.
	
	\vspace{2mm}
	
	\item If $d_n \ge \sqrt{3}\varepsilon N$, $(u_{n+1}, T_{n+1}) = (u_{n} - \Delta u, T_{n})$ and mark it by circle, go to Step 5.
	
	\vspace{2mm}
	
	\item Terminate when $u_{n+1} < 0$, otherwise go back to Step 2 with parameters $(u_{n+1}, T_{n+1})$.
\end{enumerate} 

	\vspace{2mm}
	
\textbf{Output:} The path $(u_0, T_0), \dots, (u_n, T_n), \dots $

\end{algorithm}

%... \note{Yunfeng: please add the searching algorithm here}

\begin{figure}[!h]
\centering
\includegraphics[width=4.8in,height=3.2in]{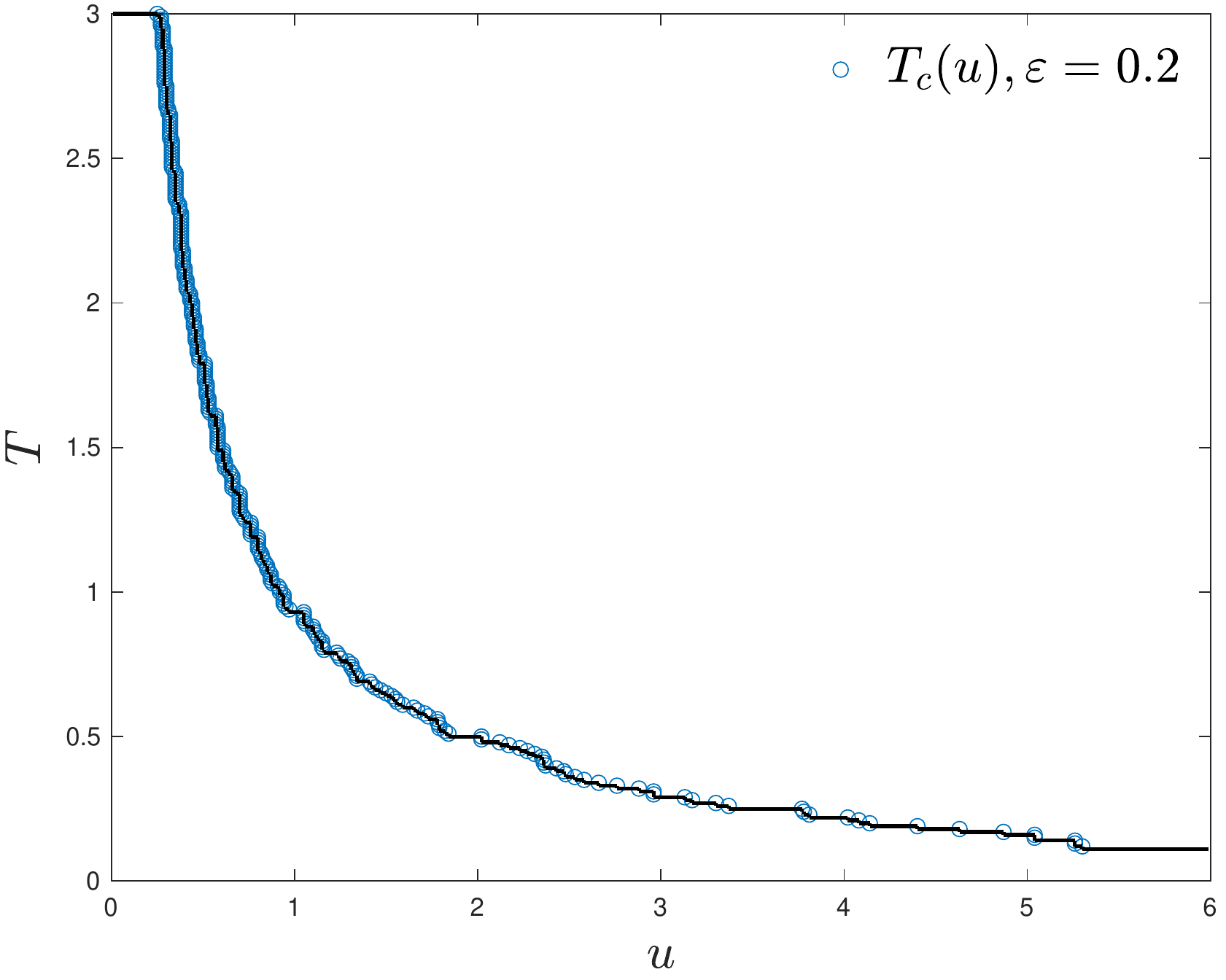}
\caption{A more precise estimation of the curve of $T_c$ via the hill climbing algorithm, $N=50, \varepsilon=0.2$.}
\label{fig_T_c}
\end{figure}

%\note{I intuitively believe $T_c(u)\asymp u^{-1}$. But if it is numerically proved wrong, then we need to re-write the following part!!}

Moreover the following linear regression on $(\log u,\log(T_c(u)))$ (marked by circle), with $N=50,\varepsilon=0.2$ indicates that, as $u\to\infty$, the lower bound of $T_c$ in Claim (iv) Theorem \ref{thm_sup_part} seems asymptotically sharp as $u\to\infty$. 

\begin{figure}[!h]
\centering
\includegraphics[width=4.8in,height=3.2in]{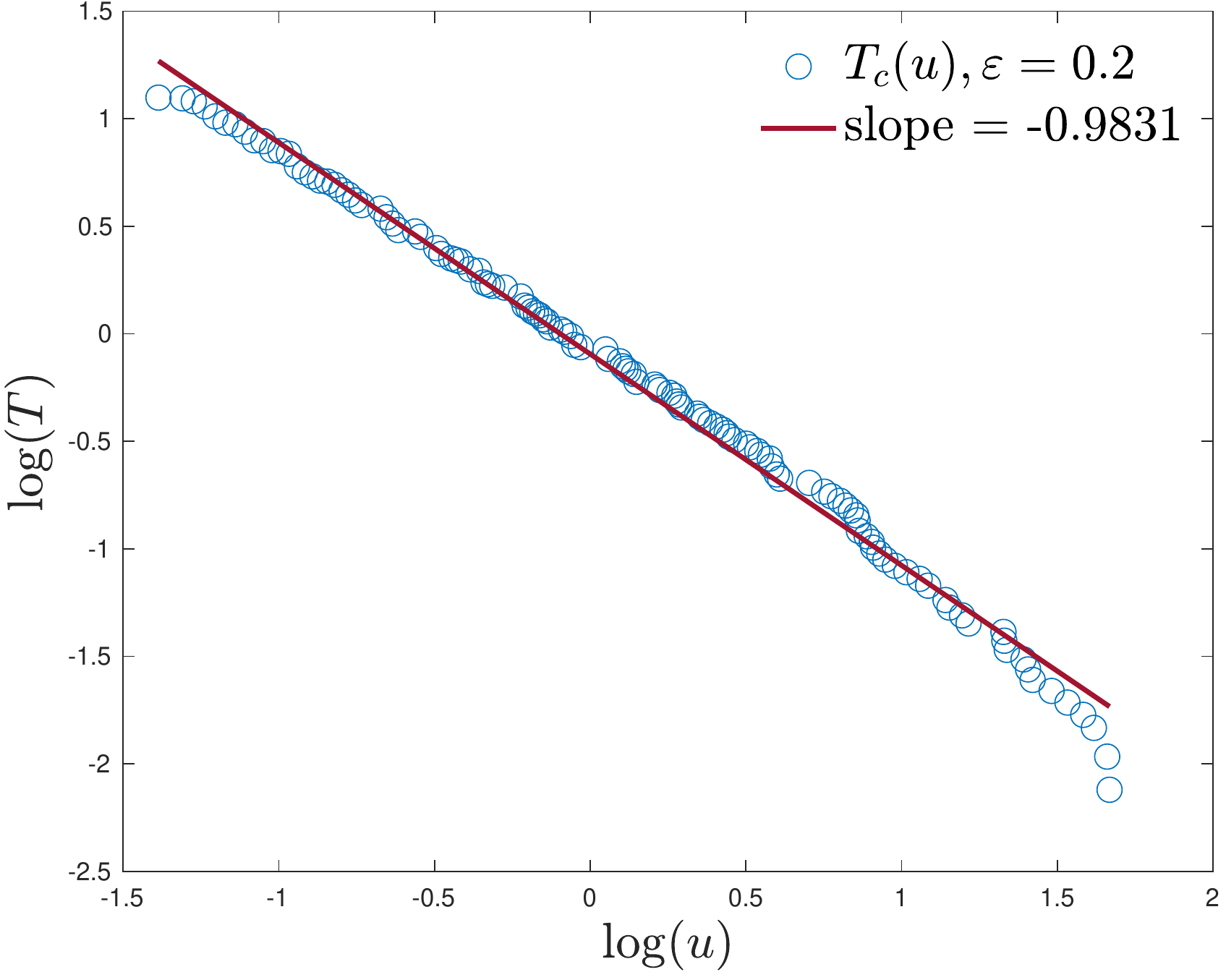}
\caption{Linear regression on $(u,T_c(u))$ in logarithm scale, $\epsilon=0.2$.The critical $T_c$ seems to be inverse proportional to $u$.}
\label{fig_reg}
\end{figure} 

\begin{conjecture}
\label{conjecture_asym}
For all $d\ge 3$, there is a constant $c_d$ such that $\lim_{u\to\infty} u T_c(u)=c_d$.
\end{conjecture}

\begin{Remark}
Simulations have also been done for different values of threshold $\epsilon$ in Algorithm \ref{hill_climbing}. All simulations strongly indicate that $T_c$ is a polynomial with respect to $u$. On the other hand, the estimated slope in linear regression under logarithm scale can be $-0.83$ when $\epsilon=0.3$, and $-0.71$ when $\epsilon=0.5$. Recalling the theoretical lower bound of $T_c$ in Claim (iv) Theorem \ref{thm_sup_part}, we believe these differences are caused by the finite size effect of $N$ expect slopes converge to $-1$, as $N\to\infty$ for all $\epsilon\in (0,1/\sqrt{3})$.
\end{Remark}

\section*{Acknowledgement}
The authors would like to thank Drs. Xinxin Chen, and Hui He, Xinyi Li, Eviatar. B. Procaccia, and Hao Wu for fruitful discussions. This work has been partially supported by Beijing Academy of Artificial Intelligence(BAAI). Y. Xiong is partially supported by the Project funded by China Postdoctoral Science Foundation (Nos. 2020TQ0011).
All the simulations performed via our own Fortran implementations run on the computing platform: 2*Intel Xeon Gold 6148 CPU  (2.40GHz, 27.5MB Cache, 40 Cores, 80 Threads) with 512GB Memory. The authors would like to sincerely thank Prof. Pingwen Zhang for allowing us accessing the computational resources of his team. This greatly accelerates our numerical simulation.

%%%%%%%%%%%%%%%%%%%%%%%%%%%%%%%%%%%%%%%%%%
%\section{Patents}
%This section is not mandatory, but may be added if there are patents resulting from the work reported in this manuscript.

%%%%%%%%%%%%%%%%%%%%%%%%%%%%%%%%%%%%%%%%%%
\vspace{6pt} 

\appendixtitles{no} % Leave argument "no" if all appendix headings stay EMPTY (then no dot is printed after "Appendix A"). If the appendix sections contain a heading then change the argument to "yes".
\appendix

\section{Proof of Theorem \ref{thm_d_large}}
\label{app_A}

In order to prove Theorem \ref{thm_d_large}, we first present a "high dimensional" version of Proposition \ref{Lawler_11_1} as follows: 

\begin{Proposition}\label{propA1}
	For $d\ge 5$, let 
	$$
	R_d=E_0[H_{d,\{0,x_1\}}\ind_{H_{d,\{0,x_1\}}<\infty}]<\infty.
	$$
	Then $\lim_{d\to\infty} R_d=0$.
\end{Proposition}

In the following arguments, we need a more precise construction of SRW on $\mathbb{Z}^d$: For $1\le i\le d$, let $\{X_n^i\}_{n=1}^{\infty}$ be an i.i.d. sequence of random variables with distribution $P(X_n^i=-1)=P(X_n^i=1)=\frac{1}{2}$. Then $\left\lbrace S_n=\sum_{k=1}^{n}X_k^i \right\rbrace_{n=0}^\infty$ forms $d$ independent copies of 1-dimension SRW's. We also define an i.i.d. sequence of random variables $\left\lbrace D_n \right\rbrace_{n=1}^\infty $ with distribution $P(D_n=j)=\frac{1}{d}$, $1\le j\le d$. Then $S_n^d=\sum\limits_{k=1}^nX_k^{D_k}e_{D_k}$ is a SRW on $\mathbb{Z}^d$. 

Before we prove Proposition \ref{propA1}, there are some preparations.

\begin{Lemma}\label{lemmaa1}
	For any $n_0>d^{1.5}$, consider a stopping time $\bar{\Gamma}:=\inf\left\lbrace n\ge 0:\exists 1\le i\le d, |S_{\tau_n^{(i)}}^{d,i}|\ge n_0^{\frac{1}{8}} \right\rbrace $, where $S_{m}^{d,i}=\sum_{k=1}^{m}X_k^i e_i$ and $\tau_n^{(i)}=\left|\left\lbrace 1\le k\le n: D_k=i\right\rbrace \right|$ for $1\le i\le d$. Then there exists $\delta>0$ such that for any sufficiently large $n_0$,\begin{equation}
		P_0\left(\bar{\Gamma}>n_0 \right)\le d*\exp(-\delta n_0^{\frac{1}{12}}).
	\end{equation}
\end{Lemma}
\begin{proof}
	By the invariance principle, we have: for any $l\ge 1$, \begin{equation}
		\lim\limits_{n_0\to \infty}P_0\left(\left| S_{ln_0^{\frac{1}{4}}}^{d,1}-S_{(l-1)n_0^{\frac{1}{4}}}^{d,1}\right|\le 2n_0^{\frac{1}{8}}\right)=P\left(B_1\le 2\right)<1,
	\end{equation}
	where $\{B_t\}_{t\ge 0}$ is a Brownian Motion starting from $0$. Therefore, there exists $\delta>0$ such that for all sufficient large $n_0$, \begin{equation}\label{a6}
		P_0\left(\left| S_{ln_0^{\frac{1}{4}}}^{d,1}-S_{(l-1)n_0^{\frac{1}{4}}}^{d,1}\right|\le 2n_0^{\frac{1}{8}}\right)<e^{-\delta}.
	\end{equation}
	
	Obviously, since $n_0>d^{1.5}$, there must exist $j_0\in \left\lbrace 1,2,...,d \right\rbrace $ such that $|\left\lbrace 1\le i\le n_0:D_i=j_0 \right\rbrace|>\frac{n_0}{d}>n_0^{\frac{1}{3}}$. By (\ref{a6}) and symmetry, we have \begin{equation}
		\begin{split}
			P_0\left(\bar{\Gamma}>n_0 \right)\le& \sum_{k=1}^{d}P_0\left(|\left\lbrace 1\le i\le n_0:D_i=j_0 \right\rbrace|>n_0^{\frac{1}{3}}, \max_{j\le n_0^{\frac{1}{3}}}|\{S_j^{d,k}\}|\le n_0^{\frac{1}{8}} \right)\\
			\le &d*P_0\left(\bigcap\limits_{l=1}^{n_0^{\frac{1}{12}}}\left\lbrace\left| S_{ln_0^{\frac{1}{4}}}^{d,1}-S_{(l-1)n_0^{\frac{1}{4}}}^{d,1}\right|\le 2n_0^{\frac{1}{8}}\right\rbrace\right)\\
			= &d*\prod_{l=1}^{n_0^{\frac{1}{12}}}P_0\left(\left| S_{ln_0^{\frac{1}{4}}}^{d,1}-S_{(l-1)n_0^{\frac{1}{4}}}^{d,1}\right|\le 2n_0^{\frac{1}{8}}\right)\le d*\exp(-\delta n_0^{\frac{1}{12}}).
		\end{split}
	\end{equation} 
\end{proof}

\begin{Lemma}\label{lemmaa2}
	There exists $c>0$ ($c$ is independent to $d$) such that for any $d\ge 20$ and any $x_0\in \mathbb{Z}^d$ such that $\exists 1\le i\le d$, $x_0^{(i)}>n_0^{\frac{1}{8}}$, \begin{equation}
		P_{x_0}\left(H_{d,\{0\}}<\infty\right)\le cn_0^{-2}.
	\end{equation}
\end{Lemma}

\begin{proof}
	Without loss of generality, we assume that $i\le \frac{d}{2}$. Then we define $I_l=\max\{1,i-9\}$, $I_r=\max\{19,i+9\}$ and $\hat{x}_0=\left(x_0^{(I_l)},x_0^{(I_l+1)},..,x_0^{(I_r)} \right)\in \mathbb{Z}^{19} $, where $||\hat{x}_0||_2\ge n_0^{\frac{1}{8}}$.
	
	We define a stopping time $\hat{H}_{d,\{0\}}=\inf\left\lbrace n\ge 0:S_{\tau_n^{(i)}}^{d,i}=0,\forall i\in [I_l,I_r]\right\rbrace $. It's easy to see that $\hat{H}_{d,\{0\}}\le H_{d,\{0\}}$. By Proposition 6.5.1 of \cite{lawler2010random}, we have \begin{equation}
		P_{x_0}\left(H_{d,\{0\}}<\infty\right)\le P_{x_0}\left(\hat{H}_{d,\{0\}}<\infty\right)=P_{\hat{x}_0}\left(H_{19,\{0\}}<\infty\right)\le c\left(n_0^{\frac{1}{8}}\right)^{-(19-2)}\le cn_0^{-2}.
	\end{equation}
\end{proof}

Now we are able to prove Proposition \ref{propA1}.

\begin{proof}[Proof of Proposition \ref{propA1}]
	Since $\left\lbrace n
	\le H_{d,\{0,x_1\}} < \infty\right\rbrace\subset \left\lbrace n
	\le H_{d,\{0\}} < \infty \right\rbrace\cup \left\lbrace n
	\le H_{d,\{x_1\}} < \infty\right\rbrace   $ for any integer $n\ge 1$, we have \begin{equation}\label{a1}
		\begin{split}
			R_d\le \sum_{n=1}^{\infty}\left[P_0\left(n
			\le H_{d,\{0\}} < \infty \right)+P_0\left(n
			\le H_{d,\{x_1\}} < \infty \right)  \right]. 
		\end{split}
	\end{equation}
	By symmetry, we have: for any $n\ge 1$,  \begin{equation}\label{a2}
		P_0\left(n+1\le H_{d,\{0\}} < \infty \right)=\frac{1}{2d}\sum_{x\in \mathbb{Z}^d, ||x||_1=1}P_x\left(n\le H_{d,\{0\}} < \infty \right)=P_{0}\left(n\le H_{d,\{x_1\}} < \infty \right).
	\end{equation}
	Combine (\ref{a1}) and (\ref{a2}), 
	\begin{equation}\label{a3}
		\begin{split}
			R_d\le& 2\sum_{n=1}^{\infty}P_0\left(n
			\le H_{d,\{0\}} < \infty \right)\\
			=&2\left[2P_0\left(2
			\le H_{d,\{0\}} < \infty \right)+\sum\limits_{n=3}^{d^{1.5}}P_0\left(n
			\le H_{d,\{0\}} < \infty \right)+\sum\limits_{n=d^{1.5}+1}^{\infty}P_0\left(n
			\le H_{d,\{0\}} < \infty \right) \right].
		\end{split}
	\end{equation}
	By the corollary of (1.14) in \cite{griffeath1983binary}, there exists $c>0$ such that for all sufficient large $d$, \begin{equation}\label{a4}
		P_0\left(H_{d,\{0\}}<\infty\right)\le \frac{1}{2d}+\frac{c}{d^2}.
	\end{equation}
	Note that $P_0\left(H_{d,\{0\}}=1\right)=0$ and $P_0\left(H_{d,\{0\}}=2\right)=\frac{1}{2d}$. By (\ref{a4}), we have\begin{equation}\label{a5}
		P_0\left(3\le H_{d,\{0\}}<\infty\right)\le \frac{c}{d^2}.
	\end{equation}
	Therefore, for the first and second term on the RHS of (\ref{a3}), we have \begin{equation}\label{a13}
		2P_0\left(2
		\le H_{d,\{0\}} < \infty \right)+\sum\limits_{n=3}^{d^{1.5}}P_0\left(n
		\le H_{d,\{0\}} < \infty \right)\le \frac{1}{d}+\frac{c}{d^{0.5}}.
	\end{equation}
	
	For the last term on the RHS of (\ref{a3}), by Lemma \ref{lemmaa1}, Lemma \ref{lemmaa2} and the strong Markov property, we have \begin{equation}\label{a14}
		\begin{split}
			\sum\limits_{n=d^{1.5}+1}^{\infty}P_0\left(n\le H_{d,\{0\}} < \infty \right)\le& \sum\limits_{n=d^{1.5}+1}^{\infty}\left[ P_0\left(\bar{\Gamma}\le n,n\le H_{d,\{0\}} < \infty\right)+P_0\left(\bar{\Gamma}> n\right)\right] \\
			\le &\sum\limits_{n=d^{1.5}+1}^{\infty}\left[cn^{-2}+d*e^{-\delta n^{\frac{1}{12}}} \right]\le c'\left( d^{-1.5}+d*e^{-\delta d^{\frac{1}{8}}}\right).
		\end{split}
	\end{equation}
	Combine (\ref{a3}), (\ref{a13}) and (\ref{a14}), \begin{equation}\label{a16}
		R_d\le 2\left[\frac{1}{d}+\frac{c}{d^{0.5}}+c'\left( d^{-1.5}+d*e^{-\delta d^{\frac{1}{8}}}\right)\right].
	\end{equation}
	Since the RHS of (\ref{a16}) converges to $0$ as $d\to \infty$, we finally get $\lim\limits_{d\to \infty}R_d=0$.
\end{proof}

Now based on Proposition \ref{propA1} and the calculations in Theorem \ref{d34}, we can get Theorem \ref{thm_d_large} easily. 

\begin{proof}[Proof of Theorem 2]
	For any $d\ge 3$, note that there exists a coupling between SRW's on $\mathbb{Z}^{d}$ and $\mathbb{Z}^{3}$: let $S_n^d=\sum\limits_{k=1}^nX_k^{D_k}e_{D_k}$ and $S_n^{3}=\sum\limits_{k=1}^nX_k^{D_k}e_{D_k}\mathbbm{1}_{1\le D_k\le 3} $. It's easy to see that $\left\lbrace H_{d,\{0,x_1\}}<\infty  \right\rbrace\subset \left\lbrace H_{3,\{0,x_1\}}<\infty \right\rbrace$. Thus, \begin{equation}\label{A1}
		Es^{(T)}_{d,\{0,x_1\}}(0)>Es_{d,\{0,x_1\}}(0)\ge Es_{3,\{0,x_1\}}(0).
	\end{equation}
	
	By Proposition \ref{propA1}, there exists $d_0$ such that for any $d>d_0$,\begin{equation}\label{A2}
		E_0^{(T)}[H_{d,\{0,x_1\}}\ind_{H_{d,\{0,x_1\}}<\infty}]<E_0[H_{d,\{0,x_1\}}\ind_{H_{d,\{0,x_1\}}<\infty}] \le Es_{3,\{0,x_1\}}(0).
	\end{equation} 
	
	By (\ref{312}), (\ref{A1}) and (\ref{A2}), for any $d>d_0$ and $0<T<\infty$, we have
	\begin{equation}
		\left[ Es_{d,\{0,x_1\}}^{(T)}(0)\left(1-P_{0}^{(T)}\left(F\right) \right)\right]'>\frac{\frac{1}{2d}}{(aT+1)^2}\left[Es_{3,\{0,x_1\}}(0)-Es_{3,\{0,x_1\}}(0)\right]=0.
	\end{equation}
	Recall that $p_{d,u}(T)=1-\exp(-4du*Es_{d,\{0,x_1\}}^{(T)}(0)\left(1-P_{0}^{(T)}\left(F\right) \right))$, then the proof is complete.
\end{proof}

%\section{}
%All appendix sections must be cited in the main text. In the appendixes, Figures, Tables, etc. should be labeled starting with `A', e.g., Figure A1, Figure A2, etc. 

%%%%%%%%%%%%%%%%%%%%%%%%%%%%%%%%%%%%%%%%%%
\reftitle{References}


%\externalbibliography{yes}
\bibliography{cite}

\begin{thebibliography}{-------}
\providecommand{\natexlab}[1]{#1}

\bibitem[Bowen(2019)]{bowen2019finitary}
Bowen, L.
\newblock Finitary random interlacements and the Gaboriau--Lyons problem.
\newblock {\em Geometric and Functional Analysis} {\bf 2019}, {\em
  29},~659--689.

\bibitem[Sznitman(2009)]{Sznitman2009Vacant}
Sznitman, A.S.
\newblock Vacant Set of Random Interlacements and Percolation.
\newblock {\em Annals of Mathematics} {\bf 2009}, {\em 171},~2039--2087.

\bibitem[Sznitman(2012)]{sznitman2012topics}
Sznitman, A.S.
\newblock {\em Topics in occupation times and Gaussian free fields}; Vol.~16,
  European Mathematical Society,  2012.

\bibitem[Teixeira(2009)]{teixeira2009interlacement}
Teixeira, A.
\newblock Interlacement percolation on transient weighted graphs.
\newblock {\em Electronic Journal of Probability} {\bf 2009}, {\em
  14},~1604--1627.

\bibitem[Procaccia \em{et~al.}(2019)Procaccia, Ye, and Zhang]{FRI_1}
Procaccia, E.B.; Ye, J.; Zhang, Y.
\newblock Percolation for the Finitary Random interlacements.
\newblock {\em arXiv preprint arXiv:1908.01954} {\bf 2019}.

\bibitem[Cai \em{et~al.}(2020)Cai, Han, Ye, and Zhang]{FRI_2}
Cai, Z.; Han, X.; Ye, J.; Zhang, Y.
\newblock On Chemical Distance and Local Uniqueness of a Sufficiently
  Supercritical Finitary Random Interlacement.
\newblock {\em arXiv preprint arXiv:2009.04044} {\bf 2020}.

\bibitem[Drewitz \em{et~al.}(2014)Drewitz, R{\'a}th, and
  Sapozhnikov]{drewitz2014chemical}
Drewitz, A.; R{\'a}th, B.; Sapozhnikov, A.
\newblock On chemical distances and shape theorems in percolation models with
  long-range correlations.
\newblock {\em Journal of Mathematical Physics} {\bf 2014}, {\em 55},~083307.

\bibitem[Liggett(2012)]{liggett2012interacting}
Liggett, T.M.
\newblock {\em Interacting particle systems}; Vol. 276, Springer Science \&
  Business Media,  2012.

\bibitem[Lawler and Limic(2010)]{lawler2010random}
Lawler, G.F.; Limic, V.
\newblock {\em Random walk: a modern introduction}; Vol. 123, Cambridge
  University Press,  2010.

\bibitem[Grimmett \em{et~al.}(2013)Grimmett, Holroyd, and
  Kozma]{grimmettpercolation}
Grimmett, G.R.; Holroyd, A.E.; Kozma, G.
\newblock Percolation of finite clusters and infinite surfaces.
\newblock {\em Mathematical Proceedings of the Cambridge Philosophical Society}
  {\bf 2013}, {\em 156},~263--279.

\bibitem[Grimmett(1999)]{MR1707339}
Grimmett, G.
\newblock {\em Percolation}, second ed.; Vol. 321, {\em Grundlehren der
  Mathematischen Wissenschaften [Fundamental Principles of Mathematical
  Sciences]}, Springer-Verlag, Berlin,  1999; pp. xiv+444.
\newblock
  doi:{\changeurlcolor{black}\href{https://doi.org/10.1007/978-3-662-03981-6}{\detokenize{10.1007/978-3-662-03981-6}}}.

\bibitem[Teixeira and Windisch(2011)]{MR2838338}
Teixeira, A.; Windisch, D.
\newblock On the fragmentation of a torus by random walk.
\newblock {\em Communications on Pure and Applied Mathematic} {\bf 2011}, {\em
  64},~1599--1646.
\newblock
  doi:{\changeurlcolor{black}\href{https://doi.org/10.1002/cpa.20382}{\detokenize{10.1002/cpa.20382}}}.

\bibitem[Griffeath \em{et~al.}(1983)Griffeath et~al.]{griffeath1983binary}
Griffeath, D.; others.
\newblock The binary contact path process.
\newblock {\em The Annals of Probability} {\bf 1983}, {\em 11},~692--705.

\end{thebibliography}

%=====================================
% References, variant B: internal bibliography
%=====================================
%\begin{thebibliography}{999}
%% Reference 1
%\bibitem[Author1(year)]{ref-journal}
%Author1, T. The title of the cited article. {\em Journal Abbreviation} {\bf 2008}, {\em 10}, 142--149.
%% Reference 2
%\bibitem[Author2(year)]{ref-book}
%Author2, L. The title of the cited contribution. In {\em The Book Title}; Editor1, F., Editor2, A., Eds.; Publishing House: City, Country, 2007; pp. 32--58.
%\end{thebibliography}

% The following MDPI journals use author-date citation: Arts, Econometrics, Economies, Genealogy, Humanities, IJFS, JRFM, Laws, Religions, Risks, Social Sciences. For those journals, please follow the formatting guidelines on http://www.mdpi.com/authors/references
% To cite two works by the same author: \citeauthor{ref-journal-1a} (\citeyear{ref-journal-1a}, \citeyear{ref-journal-1b}). This produces: Whittaker (1967, 1975)
% To cite two works by the same author with specific pages: \citeauthor{ref-journal-3a} (\citeyear{ref-journal-3a}, p. 328; \citeyear{ref-journal-3b}, p.475). This produces: Wong (1999, p. 328; 2000, p. 475)

%%%%%%%%%%%%%%%%%%%%%%%%%%%%%%%%%%%%%%%%%%
%% optional
%\sampleavailability{Samples of the compounds ...... are available from the authors.}

%% for journal Sci
%\reviewreports{\\
%Reviewer 1 comments and authors’ response\\
%Reviewer 2 comments and authors’ response\\
%Reviewer 3 comments and authors’ response
%}

%%%%%%%%%%%%%%%%%%%%%%%%%%%%%%%%%%%%%%%%%%

\end{document}